\documentclass[11pt, letterpaper]{amsart}
\usepackage{fullpage,amsfonts,amsmath,amssymb,amsthm,tikz}
\usetikzlibrary{matrix}
\usepackage[colorlinks=true,linkcolor=blue!50!black,citecolor=blue!50!black,urlcolor=blue!50!black,pdftitle={An explicit basis for the unipotent Specht modules of GLn(Fq)},pdfauthor={Scott Andrews}]{hyperref}
\let\amstext\text
\renewcommand{\text}[1]{\amstext{\upshape #1}}
\newtheorem{theorem}{Theorem}
\newtheorem{proposition}{Proposition}
\newtheorem{lemma}{Lemma}
\newtheorem{corollary}{Corollary}
\theoremstyle{definition}
\newtheorem{definition}{Definition}
\newtheorem{example}{Example}
\theoremstyle{remark}
\newtheorem{remark}{Remark}
\newcommand{\FF}{\mathbb{F}}
\newcommand{\KK}{\mathbb{K}}
\newcommand{\CC}{\mathbb{C}}
\newcommand{\ZZ}{\mathbb{Z}}
\newcommand{\GL}{\mathrm{GL}}
\newcommand{\sgn}{\mathrm{sgn}}
\newcommand{\row}{\mathrm{row}}
\newcommand{\col}{\mathrm{col}}
\newcommand{\Stab}{\mathrm{Stab}}
\newcommand{\rank}{\mathrm{rank}}
\newcommand{\rev}{\mathrm{rev}}
\newcommand{\std}{\mathrm{std}}
\newcommand{\maj}{\mathrm{maj}}
\newcommand{\Ind}{\mathrm{Ind}}
\newcommand{\Inf}{\mathrm{Inf}}
\newcommand{\St}{\mathrm{St}}

\newcommand{\domeq}{\trianglelefteq}
\newcommand{\ytab}[1]{\begin{tikzpicture}[baseline=(M.center)]
\matrix (M)[matrix of math nodes,ampersand replacement=\&,nodes={draw,minimum size=6mm,inner sep=0pt,anchor=center},column sep=-\pgflinewidth,row sep=-\pgflinewidth]{#1\\};
\end{tikzpicture}}

\begin{document}
\title{An explicit basis for the unipotent Specht modules of $\GL_n(\FF_q)$}
\author{Scott Andrews}
\address{Department of Mathematics \\ Boise State University \\ Boise, ID 83725}
\email{scottandrews@boisestate.edu}
\date{}
\begin{abstract}
For each partition $\lambda$ of $n$ we construct an explicit basis of the unipotent Specht module $S^\lambda$ of $\GL_n(\FF_q)$, defined recursively over horizontal strips. Each basis vector is $g\,x_\tau$, where $\tau$ is a semistandard tableau of shape $\lambda$, $x_\tau$ is obtained from a Steinberg-type vector by alternately applying a permutation matrix and a Fourier sum over a unitriangular group, and $g$ runs over an explicit set of coset representatives. Grouping the basis by the standardization of $\tau$ gives exactly $q^{\maj(s^T)}$ vectors for each standard tableau $s$, which is the numerical form of the Dipper--James standard basis problem. The same family is a basis over every field containing a nontrivial $p$th root of unity, including fields whose characteristic divides $|\GL_n(\FF_q)|$, and a subfamily chosen layer by layer gives a basis of the simple quotient $D^\lambda$. We also give a basis of the $U$-fixed vectors of $S^\lambda$ and computed values of $\dim D^\lambda$ for $\GL_n(2)$.
\end{abstract}
\maketitle

\section{Introduction}

For the symmetric group $S_n$, James's construction attaches to each partition $\lambda$ of $n$ a Specht module $S^\lambda$ spanned by the polytabloids $e_t$, one for each tableau $t$ of shape $\lambda$. The polytabloids $e_t$ for the standard tableaux $t$ form a basis of $S^\lambda$ over every field, so $\dim S^\lambda$ is the number of standard tableaux. Over a field of characteristic $0$ the $S^\lambda$ are the irreducible modules. In characteristic $\ell>0$ the irreducible modules are the nonzero quotients $D^\lambda=S^\lambda/(S^\lambda\cap(S^\lambda)^\perp)$ of the Specht modules by the radical of a natural bilinear form, and $D^\lambda\ne0$ for exactly the $\ell$-regular partitions $\lambda$.

The finite general linear group $G=\GL_n(\FF_q)$ has unipotent modules $S^\lambda$ indexed by the partitions $\lambda$ of $n$, the analogue of the Specht modules of $S_n$. James~\cite{James84} constructed them over every field of characteristic not dividing $q$. His method differs from the tableau method for $S_n$: each module is generated by a single explicit element attached to the initial tableau, whereas the Specht module of $S_n$ is spanned by vectors attached to all tableaux. In~\cite{AN2} the author gave a construction that is closer to the one for $S_n$. The boxes of a Young diagram are filled with linearly independent vectors of $\FF_q^n$ (an $\FF_q^n$-tableau), the group $G$ acts on these tableaux, and the Specht module $S^\lambda$ is spanned by vectors $e_T$ attached to the tableaux $T$ of shape $\lambda$ (Section~\ref{secback}). Over $\CC$, $S^\lambda$ is irreducible and affords the unipotent character $\chi^\lambda$. Over a field of characteristic $\ell\ne p$ it has a simple quotient $D^\lambda$, defined as for $S_n$ by the radical of a bilinear form. For $\GL_n(\FF_q)$ this quotient is nonzero for every partition, whereas for $S_n$ one needs $\lambda$ to be $\ell$-regular.

For $S_n$ the construction is completed by a basis theorem. For $\GL_n(\FF_q)$ it is not known how to choose, from among the vectors $e_T$, a basis of $S^\lambda$ that is described combinatorially. A related question is the standard basis problem of Dipper and James~\cite{DJ}. We have
\[
\dim S^\lambda=\sum_{t}q^{f(t)},
\]
where $t$ runs over the standard tableaux of shape $\lambda'$ and $f(t)$ is the major index of $t$~\cite{Thomas}. One would like to attach to each such $t$ a set of $q^{f(t)}$ vectors of $S^\lambda$, naturally defined and specializing as $q\to1$ to James's standard polytabloid, so that together they form a basis. Dipper and James made conjectures about the form of the answer and verified them for the two-part partitions $(n-m,m)$ in James's labeling (two-column shapes in ours) with $m\le11$, using computations of Brandt. The two-part case was proved by Brandt, Dipper, James, and Lyle~\cite{BDJL}, and Guo~\cite{Guo,DG} gave a representation-theoretic construction of the standard basis in that case. For general $\lambda$ no basis of this kind is known to us.

\medskip
\noindent\textbf{Results.} We construct a basis of $S^\lambda$ by recursion on the horizontal strips of $\lambda$. Restrict $S^\lambda$ to the mirabolic subgroup, the matrices of $G$ with last row $(0,\dots,0,1)$ (Example~\ref{exmirab}). A module $X$ of this subgroup has layers $\Omega_{Y_k}X$, $1\le k\le n$, which are determined by the action of a unitriangular group on the last $k$ columns. Each layer is a module for a smaller group $\GL_{n-k}$, and $X$ is the direct sum of the translates of the layers by an explicit set of coset representatives (Lemma~\ref{lemlayers}). For $X=S^\lambda$, Harish-Chandra theory identifies the $k$th layer with the direct sum of the Specht modules $S^\mu$ of $\GL_{n-k}$ over the horizontal $k$-strips $\lambda/\mu$ (Proposition~\ref{proplayer}).

Suppose bases of the $S^\mu$ are known. For each $\mu$ we choose a single vector in the $k$th layer whose component in the summand $S^\mu$ is a generator of $S^\mu$, and whose components in the other summands $S^\nu$ vanish unless $\nu\triangleright\mu$. Applying the known basis of $S^\mu$ to this vector, for all $\mu$, gives a basis of the layer, and repeating the construction gives a basis of $S^\lambda$ defined by recursion on $n$.

Explicitly (Definition~\ref{defB} and Theorem~\ref{thmmain}), let $x_{T_\lambda}\in S^\lambda$ be the analogue of Steinberg's vector $\sum_w\varepsilon(w)\,w\,m_T$ attached to the row-reading tableau $T_\lambda$. The basis of $S^\lambda$ consists of the vectors $\beta\,x_{T_\lambda}$ with $\beta\in B(\lambda)\subset\CC[G]$, where $B(\lambda)$ is built from the following data:
\begin{itemize}
\item a chain $\lambda=\nu^0\supset\nu^1\supset\dots\supset\nu^r=\emptyset$ of partitions in which each $\nu^j/\nu^{j+1}$ is a horizontal strip, that is, a semistandard tableau $\tau$ of shape $\lambda$;
\item for each step, a permutation matrix $W$ that moves the boxes of the strip to the end, and the Fourier sum $\Omega_{Y_k}$ over a unitriangular group, a scalar multiple of the projection onto a character eigenspace (for a diagram with $m$ boxes and a strip with $k$ boxes, a sum of $q^{k(m-k)+\binom k2}$ group elements);
\item for each step, one element of a set $R_k$ of coset representatives (Section~\ref{subnotation}).
\end{itemize}
Thus the basis is indexed by a semistandard tableau $\tau$ of shape $\lambda$ together with a choice of coset representatives, and the number of choices depends only on the content of $\tau$ (Corollary~\ref{corclosed}). We prove the following.
\begin{itemize}
\item The vectors form a basis of $S^\lambda$ over $\CC$ (Theorem~\ref{thmmain}). The proof also gives the constant relating the vectors $x_T$ and $e_T$ (Corollary~\ref{corkappa}) and a triangularity theorem for the layers (Theorem~\ref{thmT}).
\item Grouping the basis by the standardization $s$ of $\tau$ gives exactly $q^{\maj(s^T)}$ vectors for each standard tableau $s$ (Proposition~\ref{propstd}), which is the numerical form of the standard basis problem. We do not claim that the vectors attached to $s$ are ``natural'' in the sense of Dipper and James, who report evidence that no natural choice exists for these numbers (Section~\ref{secstd}).
\item Let $\KK$ be a field containing a nontrivial $p$th root of unity, so that its characteristic $\ell$ is $0$ or different from $p$ ($\ell$ may divide $|G|$). The same family is linearly independent over $\KK$ (Theorem~\ref{thmindep}) and is a basis of $S^\lambda_\KK$, which is James's Specht module (Theorem~\ref{thmP}).
\item The images of a subfamily of the recursive basis, chosen layer by layer, form a basis of the simple module $D^\lambda$ (Proposition~\ref{propD}). The subfamily is determined by the radical of the bilinear form, as in the case of the symmetric group.
\item The vectors $x_\tau$ attached to the standard chains $\tau$ of $\lambda$ are fixed by the Borel subgroup and form a basis of the space of $U$-fixed vectors of $S^\lambda_\KK$, whose dimension is the number of standard tableaux (Proposition~\ref{propU}). Their Gram matrix has rank $\dim(D^\lambda)^U$ (Proposition~\ref{propGramU}). Section~\ref{subsym} explains how this relates $D^\lambda$ to the Hecke algebra $H_q(S_n)$ and to the symmetric group, using the work of Dipper and James.
\item Section~\ref{seccomp} gives computed values of $\dim D^\lambda$ and $\dim(D^\lambda)^U$ for $\GL_n(2)$ with $n\le5$, and compares them with James's decomposition matrices for $n\le4$ and with the block theory of the Hecke algebra. They agree, except for one entry of James's matrix that appears to be misprinted.
\end{itemize}

\medskip
\noindent\textbf{Relation to James's book.} The layer decomposition (Lemma~\ref{lemlayers}) is~\cite[Thm.~9.11]{James84}, and the identification of the layers of $S^\lambda$ (Proposition~\ref{proplayer}, for which we give a different proof) corresponds to James's filtration of the layers by smaller Specht modules~\cite[Thm.~16.9]{James84}, which splits in characteristic $0$. His Branching Theorem~\cite[Thm.~16.10]{James84} states the resulting restriction to $G_{n-1}$. Along the same layers, \cite[Ch.~12 (p.~73), Thm.~16.9]{James84} constructs linearly independent elements of $S^\lambda$ that form a basis by his Kernel Intersection Theorem. His family is described through its leading terms (his (12.11)), and the Kernel Intersection Theorem gives the upper bound on the dimension. Our family consists of explicit vectors whose independence is proved directly (Theorem~\ref{thmindep}).

\medskip
\noindent\textbf{Organization.} Section~\ref{secback} recalls the Specht modules of~\cite{AN2}, the facts about unipotent characters that we need, and the notation for the basis. Section~\ref{secbasis} defines the basis and states the main theorem, with examples. Sections~\ref{secx}--\ref{secdiag} contain the proof: the vectors $x_T$, the layer decomposition, the layers of $S^\lambda$, the inner products with $x_T$, and the three properties (nonvanishing, triangularity, the diagonal component) of the vectors $\Omega_{Y_k}x_{T'}$ that make the induction work. Section~\ref{secproof} completes the proof, and Section~\ref{secstd} groups the basis by standard tableaux and counts the vectors in each group. Section~\ref{secmodular} extends the results to positive characteristic and to the simple modules $D^\lambda$, and gives a basis of the $U$-fixed vectors. Section~\ref{seccomp} records computed values of $\dim D^\lambda$, and Section~\ref{secfurther} lists open questions and limitations.

\section{Background}\label{secback}

\subsection{Partitions and tableaux}

Let $n$ be a positive integer. A \emph{partition} of $n$ is a sequence $\lambda=(\lambda_1,\lambda_2,\dots,\lambda_k)$ of positive integers with $\lambda_1\ge\lambda_2\ge\dots\ge\lambda_k$ and $\lambda_1+\lambda_2+\dots+\lambda_k=n$. We write $\lambda\vdash n$ to indicate that $\lambda$ is a partition of $n$.

There is a partial order on the set of partitions of $n$, called the \emph{dominance order}, with $\lambda\trianglerighteq\mu$ if and only if
\[
\sum_{i=1}^{j}\lambda_i\ge\sum_{i=1}^{j}\mu_i
\]
for all $j$ (setting $\lambda_i=0$ if $\lambda$ has fewer than $i$ parts).

To each partition $\lambda$ we associate a \emph{Young diagram}, a left-justified array of boxes in which the number of boxes in row $i$ is $\lambda_i$. A box of the diagram lies in row $r$ and column $c$.

The \emph{conjugate} of a partition $\lambda$, denoted $\lambda'$, is the partition defined by $\lambda'_i=|\{j\mid\lambda_j\ge i\}|$. The Young diagram of $\lambda'$ is obtained from that of $\lambda$ by reflection across the diagonal.

\begin{example}\label{exdiag}
For $\lambda=(4,2,1)$ the Young diagram of $\lambda$ and that of the conjugate partition $\lambda'=(3,2,1,1)$ are
\[
\lambda\colon\ \ytab{{}\&{}\&{}\&{}\\{}\&{}\\{}}\qquad\qquad\lambda'\colon\ \ytab{{}\&{}\&{}\\{}\&{}\\{}\\{}}
\]
\end{example}

If $\lambda$ is a partition of $n$, a \emph{tableau} of shape $\lambda$ is a filling of the Young diagram of $\lambda$ with positive integers. A tableau is \emph{standard} if it contains each of the integers $1$ through $n$ exactly once and its entries increase along rows and down columns, and \emph{semistandard} if its entries weakly increase along rows and strictly increase down columns.

If $\mu\subseteq\nu$ are partitions (that is, $\mu_i\le\nu_i$ for all $i$, so that the diagram of $\mu$ lies inside that of $\nu$), the \emph{skew diagram} $\nu/\mu$ is the set of boxes of $\nu$ that are not in $\mu$. It is a \emph{horizontal $k$-strip} ($k=|\nu|-|\mu|$) if it has $k$ boxes, no two in the same column; equivalently $\nu_{i+1}\le\mu_i\le\nu_i$ for all $i$.

A chain $\lambda=\nu^0\supset\nu^1\supset\dots\supset\nu^r=\emptyset$ in which every $\nu^j/\nu^{j+1}$ is a nonempty horizontal strip can be identified with the semistandard tableau $\tau$ of shape $\lambda$ that has the entry $r-j$ in the boxes of $\nu^j/\nu^{j+1}$. The entries of such a tableau are exactly $1,\dots,r$, each occurring at least once, and the \emph{content} of $\tau$ is the composition $\alpha=(\alpha_1,\dots,\alpha_r)$ of $|\lambda|$ in which $\alpha_i$ is the number of entries equal to $i$ (so $\nu^j/\nu^{j+1}$ has $\alpha_{r-j}$ boxes). We will only use semistandard tableaux of this kind, so from now on a semistandard tableau always has entries exactly $1,\dots,r$ for some $r$; in particular, its entries are positive integers and no number is skipped.

\begin{example}\label{exstrip}
Let $\lambda=(4,3,1)$ and consider the chain
\[
\lambda=\nu^0\supset\nu^1\supset\nu^2\supset\nu^3\supset\nu^4=\emptyset,\qquad \nu^1=(3,2),\ \nu^2=(2,1),\ \nu^3=(2).
\]
Each $\nu^j/\nu^{j+1}$ is a horizontal strip, of sizes $3,2,1,2$. In the picture each diagram is $\nu^j$, with the strip $\nu^j/\nu^{j+1}$ shaded:
\[
\ytab{{}\&{}\&{}\&|[fill=gray!35]|\\{}\&{}\&|[fill=gray!35]|\\|[fill=gray!35]|}\ \supset\ \ytab{{}\&{}\&|[fill=gray!35]|\\{}\&|[fill=gray!35]|}\ \supset\ \ytab{{}\&{}\\|[fill=gray!35]|}\ \supset\ \ytab{|[fill=gray!35]|\&|[fill=gray!35]|}\ \supset\ \emptyset.
\]
The corresponding semistandard tableau $\tau$ has the entry $4-j$ in the boxes of $\nu^j/\nu^{j+1}$:
\[
\tau=\ytab{1\&1\&3\&4\\2\&3\&4\\4}.
\]
\end{example}

The \emph{standardization} $\std(\tau)$ of a semistandard tableau $\tau$ of shape $\lambda\vdash n$ is the standard tableau obtained by replacing the entries of $\tau$ with the integers $1$ through $n$, starting with the smallest entries of $\tau$ and, among equal entries, moving from left to right. For the tableau $\tau$ above, $\std(\tau)$ has rows $1\,2\,5\,8$, $3\,4\,7$, and $6$:
\[
\std(\tau)=\ytab{1\&2\&5\&8\\3\&4\&7\\6}.
\]
The Kostka number $K_{\nu\lambda}$ is the number of semistandard tableaux of shape $\nu$ with $\lambda_i$ entries equal to $i$.

\subsection{The Specht modules of \texorpdfstring{$\GL_n(\FF_q)$}{GLn(Fq)}}\label{subspecht}

Let $q$ be a power of the prime $p$. The field $\FF_q$ is fixed throughout, so we write $\GL_m$ for $\GL_m(\FF_q)$ and put $G=\GL_n$. Let $\KK$ be a field that contains a nontrivial $p$th root of unity (so $\mathrm{char}\,\KK\ne p$), and fix a nontrivial homomorphism $\theta:(\FF_q,+)\to\KK^\times$. We use the construction of~\cite{AN2}. An \emph{$\FF_q^n$-tableau} of shape $\lambda\vdash n$ is a filling of the boxes of the Young diagram of $\lambda$ with $n$ linearly independent vectors of $\FF_q^n$. Let $\mathcal T_\lambda$ denote the set of $\FF_q^n$-tableaux of shape $\lambda$. The group $G$ acts on $\mathcal T_\lambda$ by left multiplication of the entries (written as column vectors), and this action is simply transitive. Let $\KK[\mathcal T_\lambda]$ be the permutation module with basis $\mathcal T_\lambda$. For $T\in\mathcal T_\lambda$ put
\begin{align*}
U(T)&=\{g\in G:\ gv-v\in\mathrm{span}\{v'\ \text{strictly left of}\ v\ \text{in}\ T\}\ \text{for every entry}\ v\},\\
P(T)&=\{g\in G:\ gv\in\mathrm{span}\{v'\ \text{weakly right of}\ v\ \text{in}\ T\}\ \text{for every entry}\ v\}.
\end{align*}
The group $U(T)$ is unipotent and $P(T)$ is a parabolic subgroup. If $X(T)$ is the set of pairs $(v,v')$ of entries such that $v$ lies in the box directly to the left of that of $v'$, then
\[
\psi_T(u)=\theta\Bigl(\sum_{(v,v')\in X(T)}(\text{the coefficient of $v$ in $uv'$})\Bigr)
\]
is a linear character of $U(T)$ (here ``the coefficient of $v$'' is taken in the basis formed by the entries of $T$). Define
\[
k_T=\sum_{u\in U(T)}\psi_T(u^{-1})\,u\in\KK[G].
\]
We also define two elements $m_T$ and $e_T$ of $\KK[\mathcal T_\lambda]$:
\[
m_T=\sum_{p\in P(T)}pT\qquad\text{and}\qquad e_T=k_T\,m_T=\sum_{u\in U(T)}\psi_T(u^{-1})\,m_{uT}.
\]
For $\lambda\vdash n$, let
\[
M^\lambda_\KK=\KK\text{-span}\{m_T: T\in\mathcal T_\lambda\}\qquad\text{and}\qquad S^\lambda_\KK=\KK\text{-span}\{e_T: T\in\mathcal T_\lambda\}.
\]
We omit the subscript when $\KK=\CC$ and write $M^\lambda$ and $S^\lambda$. These are $G$-modules, since $g\,m_T=m_{gT}$ and $g\,e_T=e_{gT}$. The groups $P(T)$ and $U(T)$ play the role of the row and column stabilizers for $S_n$, and $\psi_T$ is the analogue of the sign character.

\begin{example}\label{extab}
Let $\lambda=(2,2)$, $n=4$, and let $T$ be the tableau with rows $(v_1,v_2)$ and $(v_3,v_4)$, where $v_1,\dots,v_4$ is the standard basis:
\[
T=\ytab{v_1\&v_2\\v_3\&v_4}.
\]
Then
\[
U(T)=\left\{\begin{pmatrix}1&a&0&c\\0&1&0&0\\0&b&1&d\\0&0&0&1\end{pmatrix}\right\},\qquad
P(T)=\left\{\begin{pmatrix}*&0&*&0\\ *&*&*&*\\ *&0&*&0\\ *&*&*&*\end{pmatrix}\right\},
\]
$X(T)=\{(v_1,v_2),(v_3,v_4)\}$ and $\psi_T(u)=\theta(a+d)$.
\end{example}

The module $M^\lambda_\KK$ is the permutation module of $G$ on the \emph{flags} of type determined by $\lambda$ (see Section~\ref{subflags}); $S^\lambda_\KK$ is a submodule. The following facts are proved in~\cite{AN2} (Lemmas~1, 4, 5, and~7, and Corollary~1 of the published version) for every $\KK$ as above. We use them with $\KK=\CC$ in Sections~\ref{secx}--\ref{secproof} and for general $\KK$ in Section~\ref{secmodular}.
\begin{itemize}
\item[(A1)] $g\,m_T=m_{gT}$, $g\,e_T=e_{gT}$, and $e_{uT}=\psi_T(u)e_T$ for $u\in U(T)$; moreover $e_T\ne0$.
\item[(A2)] If $T,T'$ are $\FF_q^n$-tableaux of the same shape, then $k_T\,m_{T'}\in\KK e_T$.
\item[(A3)] If $\mathrm{char}\,\KK\nmid|G|$ (for instance $\KK=\CC$), $S^\lambda_\KK$ is irreducible.
\end{itemize}
Since $G$ acts simply transitively on tableaux of a given shape, $S^\lambda_\KK=\KK[G]\,e_T$ for any $T\in\mathcal T_\lambda$.

\subsection{Unipotent characters}\label{subunip}

The irreducible complex characters of $\GL_n$ were determined by Green~\cite{G}. We follow the labeling of Green and Macdonald~\cite{G,Mac}; James's labeling is the transpose. The \emph{unipotent} characters are the irreducible constituents of $\Ind_B^G1$ ($B$ a Borel subgroup). They are the characters $\chi^\nu$, $\nu\vdash n$, with $\chi^{(1^n)}$ trivial and $\chi^{(n)}$ the Steinberg character $\St$. By~\cite{AN2}, $S^\lambda$ affords $\chi^\lambda$. We use two classical facts.

\begin{itemize}
\item[(H1)] Harish-Chandra induction corresponds to the multiplication of symmetric functions ($\chi^\nu\mapsto s_\nu$, $\St_{\GL_n}\mapsto h_n$, $1_{\GL_n}\mapsto e_n$)~\cite{G,Z}. Consequently, for a parabolic subgroup $P_\mu\subset\GL_n$ of type $\mu$ (block sizes $\mu_1,\mu_2,\dots$), \[
\Ind_{P_\mu}^G1=\sum_{\nu'\trianglerighteq\mu}K_{\nu'\mu}\chi^\nu,\qquad \Ind_{P_\mu}^G(\St_{\mu_1}\otimes\St_{\mu_2}\otimes\cdots)=\sum_{\nu\trianglerighteq\mu}K_{\nu\mu}\chi^\nu,
\]
and the Jacquet module (Harish-Chandra restriction) of a module affording $\chi^\lambda$ to a Levi subgroup $\GL_{n-k}\times\GL_k$ affords $\sum_{\mu,\rho}c^\lambda_{\mu\rho}\,\chi^\mu\otimes\chi^\rho$, with $c^\lambda_{\mu\rho}$ the Littlewood--Richardson coefficients.
\item[(H2)] (Zelevinsky~\cite{Z}, 9.4 and 9.7) Among the unipotent modules of $\GL_k$, only the Steinberg module has a nonzero $(U_k,\psi)$-eigenspace, where $U_k$ is the group of upper unitriangular matrices and $\psi(u)=\theta(\sum_iu_{i,i+1})$ a generic character; this eigenspace is one dimensional. Equivalently, the Steinberg character is the only unipotent character that occurs in the Gelfand--Graev character $\Ind_{U_k}^{\GL_k}\psi$, and it occurs with multiplicity one.
\end{itemize}
Since $c^\lambda_{\mu,(k)}=1$ if $\lambda/\mu$ is a horizontal $k$-strip and $0$ otherwise, (H1) and (H2) say that the $(U_k,\psi)$-eigenspace of the Jacquet module of $S^\lambda$ for $\GL_{n-k}\times\GL_k$ is $\bigoplus_{\lambda/\mu\ \text{horizontal $k$-strip}}S^\mu$ as a $\GL_{n-k}$-module (see Proposition~\ref{proplayer}).

\subsection{Coordinate tableaux, flags, and the vectors \texorpdfstring{$x_T$}{xT}}\label{subflags}

Let $v_1,\dots,v_n$ be the standard basis of $\FF_q^n$ and, for $g\in\GL_n$, write $g_{ab}$ for the coefficient of $v_a$ in $gv_b$ (in other words, $g_{ab}$ is the $(a,b)$ entry of the matrix of $g$ in this basis). An $\FF_q^n$-tableau is a \emph{coordinate tableau} if its entries are $v_1,\dots,v_n$ in some order. We call $a$ the \emph{label} of $v_a$, and for a coordinate tableau $T$ we write $\row(a),\col(a)$ for the row and column of the box of $T$ that contains $v_a$. The \emph{row-reading tableau} $T_\lambda$ is the coordinate tableau with $v_1,\dots,v_{\lambda_1}$ in the first row, the next $\lambda_2$ vectors in the second row, and so on. For a coordinate tableau $T$, the groups $P(T)$ and $U(T)$ of Section~\ref{subspecht} are described by sets of positions:
\begin{align*}
P(T)&=\{g:\ g_{ab}\ne0\Rightarrow\col(a)\ge\col(b)\},\\
U(T)&=\{u:\ u_{aa}=1,\ u_{ab}\ne0\ (a\ne b)\Rightarrow\col(a)<\col(b)\}.
\end{align*}
We also introduce two ``row'' groups for an arbitrary $\FF_q^n$-tableau $T$:
\begin{align*}
N(T)&=\{g\in G:\ gv-v\in\mathrm{span}\{v'\ \text{strictly above}\ v\ \text{in}\ T\}\ \text{for every entry}\ v\},\\
R(T)&=\{g\in G:\ g\ \text{permutes the entries of}\ T\ \text{and}\ gv\ \text{lies in the row of}\ v\ \text{for every entry}\ v\}.
\end{align*}
For a coordinate tableau $T$ these are
\begin{align*}
N(T)&=\{n:\ n_{aa}=1,\ n_{ab}\ne0\ (a\neq b)\Rightarrow\row(a)<\row(b)\},\\
R(T)&=\{w:\ w\ \text{a permutation matrix},\ w_{ab}\ne0\Rightarrow\row(a)=\row(b)\}.
\end{align*}
Thus $R(T)$ is the group of permutations of the entries of $T$ that preserve every row of $T$ (it is isomorphic to $S_{\lambda_1}\times S_{\lambda_2}\times\cdots$), and we write $\sgn(r)$ for the sign of the permutation that $r\in R(T)$ induces on the entries.

Let $T$ be an $\FF_q^n$-tableau of shape $\lambda$. For $c\ge1$ put
\[
W_c(T)=\mathrm{span}\{v:\ v\ \text{lies in a column}\ \ge c\ \text{of}\ T\},
\]
so $W_1(T)=\FF_q^n$, $W_c(T)=0$ for $c>\lambda_1$, and $\dim W_c(T)=\lambda'_c+\lambda'_{c+1}+\cdots$. The \emph{flag} of $T$ is the chain
\[
F_T:\qquad \FF_q^n=W_1(T)\supsetneq W_2(T)\supsetneq\cdots\supsetneq W_{\lambda_1}(T)\supsetneq W_{\lambda_1+1}(T)=0.
\]
(We write $W_c$ for $W_c(T)$ when $T$ is clear.) Let $X_\lambda$ be the set of all flags $\FF_q^n=W_1\supsetneq\cdots\supsetneq W_{\lambda_1+1}=0$ with $\dim W_c/W_{c+1}=\lambda'_c$ (flags of type $\lambda'$). In this subsection $M^\lambda=M^\lambda_\KK$ for a field $\KK$ as in Section~\ref{subspecht}, except where we say otherwise. These flags have the following properties.
\begin{enumerate}
\item[(a)] $F_T\in X_\lambda$ and $gF_T=F_{gT}$ for $g\in G$, since $gW_c(T)=W_c(gT)$. Every flag in $X_\lambda$ is $F_T$ for some $T$: choose, for each $c$, $\lambda'_c$ vectors that span a complement of $W_{c+1}$ in $W_c$, and place them in column $c$. Hence $G$ acts transitively on $X_\lambda$.
\item[(b)] $P(T)=\mathrm{Stab}_G(F_T)$. Indeed, $gW_c\subseteq W_c$ for all $c$ says exactly that $gv$ lies in the span of the entries in columns $\ge\col(v)$, for each entry $v$ of $T$. In particular $X_\lambda\cong G/P(T)$, and $P(T)$ is a parabolic subgroup of type $\lambda'$.
\item[(c)] $m_{T'}=m_{T''}$ if and only if $F_{T'}=F_{T''}$ (by (b)). So $m_T$ depends only on the flag $F_T$, and $m_T\mapsto[F_T]$ identifies $M^\lambda$ with the permutation module $\KK[X_\lambda]\cong\Ind_{P(T)}^G1$, where $[F]$ is the basis vector of the flag $F$. From now on we identify $M^\lambda$ with $\KK[X_\lambda]$ in this way. For $x\in M^\lambda$ and a flag $F$, we write $x(F)$ for the coefficient of $[F]$ in $x$, so that $x=\sum_Fx(F)[F]$.
\item[(d)] Over $\CC$ we use the Hermitian form $\langle\ ,\ \rangle$ on $M^\lambda$ for which the flags are orthonormal. In other words, $\langle x,y\rangle=\sum_Fx(F)\overline{y(F)}$. It is $G$-invariant, since $G$ permutes the flags. (For general $\KK$ we use the symmetric bilinear form $[\,,\,]$ with the flags orthonormal; see Section~\ref{subindep}.)
\item[(e)] Over $\CC$, since $P(T)$ is a parabolic subgroup of type $\lambda'$, (H1) gives that $M^\lambda$ has character $\sum_{\nu'\trianglerighteq\lambda'}K_{\nu'\lambda'}\chi^\nu$. Hence $M^\lambda\cong\bigoplus_{\nu\domeq\lambda}(S^\nu)^{\oplus K_{\nu'\lambda'}}$; in particular $S^\lambda$ occurs exactly once, and every other constituent is $S^\nu$ with $\nu\lhd\lambda$.
\item[(f)] In this language
\[
e_T=\sum_{u\in U(T)}\psi_T(u)^{-1}\,[u\,F_T],
\]
and the flags $uF_T$ are distinct for distinct $u\in U(T)$, since $U(T)\cap P(T)=1$.
\end{enumerate}
Finally define
\[
x_T=\sum_{r\in R(T)}\sgn(r)\,r\sum_{n\in N(T)}n\,m_T\ \in M^\lambda
\]
for every $\FF_q^n$-tableau $T$ of shape $\lambda$. Since $N(gT)=gN(T)g^{-1}$, $R(gT)=gR(T)g^{-1}$, and $g\,m_T=m_{gT}$, we have $g\,x_T=x_{gT}$. For $\lambda=(n)$ we have $N(T)=1$, and $x_T=\sum_{w\in R(T)}\varepsilon(w)\,w\,m_T$ is Steinberg's vector~\cite{St} for the basis formed by the entries of $T$. The vector $x_T$ is not one of the $e_T$, but it generates the same module (Theorem~\ref{thmx} below), and it is the starting point of our recursion.

\begin{example}\label{exxT}
For $\lambda=(2)$ and $T=T_\lambda=\ytab{v_1\&v_2}$ we have $W_2(T)=\mathrm{span}\{v_2\}$, so $F_T$ is the flag $\FF_q^2\supsetneq\mathrm{span}\{v_2\}\supsetneq0$. Then $P(T)=\mathrm{Stab}_G(F_T)$ is the group of lower triangular matrices, $N(T)=1$, and $R(T)=\{1,w\}$, where $w$ is the permutation matrix that swaps $v_1$ and $v_2$. Since $m_T=[F_T]$ and $w\,m_T=[wF_T]$, we get $x_T=[F_T]-[wF_T]$. The set $X_{(2)}$ consists of the flags $\FF_q^2\supsetneq W\supsetneq0$ with $W$ one dimensional, so there are $q+1$ of them and $M^{(2)}$ has dimension $q+1$. The Steinberg module $S^{(2)}$ is the $q$-dimensional space of combinations of flags whose coefficients sum to $0$.
\end{example}

\subsection{Notation for the basis}\label{subnotation}

We use the conventions that $\GL_0$ is the trivial group, $S^\emptyset$ is one dimensional with generator $x_{T_\emptyset}=1$, empty products are $1$, and $x^+=\max(x,0)$ for a real number $x$. We regard $\GL_m\subset\GL_n$ as the top-left block ($m\le n$), so every matrix and group defined below lies in $\GL_m$ and acts on $v_1,\dots,v_m$ and fixes $v_{m+1},\dots,v_n$. Fix $m\ge1$ and $1\le k\le m$. When $m$ is clear from the context we drop it from the notation below, writing $Y_k$, $\psi_k$, and $R_k$ for $Y^{(m)}_k$, $\psi^{(m)}_k$, and $R^{(m)}_k$; we keep the superscript only where several values of $m$ occur together.

\begin{itemize}
\item $Y^{(m)}_k\subset\GL_m$ is the group of upper unitriangular $m\times m$ matrices all of whose off-diagonal entries lie in the last $k$ columns; $\psi^{(m)}_k$ is its character $\psi^{(m)}_k(y)=\theta\bigl(\sum_{i=m-k+2}^{m}y_{i-1,i}\bigr)$ (the sum of the entries $y_{i-1,i}$ just above the diagonal that lie in the last $k-1$ columns; the sum is empty for $k=1$), and
\[
\Omega_{Y_k}=\sum_{y\in Y_k}\psi_k(y)^{-1}\,y\ \in\KK[\GL_m].
\]
For a module $X$ of $Y_k$, the \emph{$(Y_k,\psi_k)$-eigenspace} of $X$ is $\{x\in X:\ yx=\psi_k(y)x\ \text{for all}\ y\in Y_k\}$. We have the following:
\begin{enumerate}
\item[(i)] $\Omega_{Y_k}^2=|Y_k|\,\Omega_{Y_k}$, so $\Omega_{Y_k}$ is a Fourier sum, not an idempotent.
\item[(ii)] $y\,\Omega_{Y_k}=\psi_k(y)\,\Omega_{Y_k}$ for $y\in Y_k$. Hence $\Omega_{Y_k}x$ lies in the eigenspace for every $x\in X$.
\item[(iii)] If $x$ lies in the eigenspace, then $\Omega_{Y_k}x=|Y_k|\,x$.
\item[(iv)] $|Y_k|$ is a power of $p$, hence invertible in $\KK$ (as $\mathrm{char}\,\KK\ne p$). So $|Y_k|^{-1}\Omega_{Y_k}$ is the idempotent projecting $X$ onto the eigenspace, and the image $\Omega_{Y_k}X$ is exactly the $(Y_k,\psi_k)$-eigenspace of $X$.
\item[(v)] Over $\CC$, $\Omega_{Y_k}$ is self-adjoint for the Hermitian form of Section~\ref{subflags}(d), because the adjoint of $y$ is $y^{-1}$ and $\overline{\psi_k(y)}=\psi_k(y)^{-1}$.
\end{enumerate}
\item For each $i\ge1$ and each nonzero row vector $w\in\FF_q^{m-i}$, choose a matrix $h_w\in\GL_{m-i}\subset\GL_m$ such that the last row of $h_w^{-1}$ is $w$. Define
$R^{(m)}_k=\{h_1h_2\cdots h_{k-1}\}$, where each $h_i$ runs over the matrices $h_w$ with $w\in\FF_q^{m-i}\setminus0$. Thus $|R_k|=\prod_{j=m-k+1}^{m-1}(q^j-1)$, and $R_1=\{1\}$. The set $R_k$ depends on the choice of the $h_w$, and different choices give different bases of $S^\lambda$ in Theorem~\ref{thmmain}. The theorem holds for every choice, since the proof of Lemma~\ref{lemlayers} uses only the property of $h_w$ stated here.
\item For a partition $\nu\vdash m$ and a horizontal $k$-strip $\nu/\mu$, $W(\nu,\mu)\in S_m\subset\GL_m$ is the permutation matrix defined as follows. Let $T_\nu$ be the row-reading tableau of $\nu$ and let $T'(\nu,\mu)$ be the coordinate tableau of shape $\nu$ that agrees with $T_\mu$ on the boxes of $\mu$ and has the labels $m-k+1,\dots,m$ in the boxes of $\nu/\mu$, in increasing order from left to right (so the smallest of these labels is in the lowest row). Then $W(\nu,\mu)v_a=v_{a'}$, where $a'$ is the label of $T'(\nu,\mu)$ in the box that carries $a$ in $T_\nu$. Thus $W(\nu,\mu)T_\nu=T'(\nu,\mu)$ as tableaux.
\end{itemize}
The labels $m-k+1,\dots,m$ of the strip are called the \emph{strip letters}.

\begin{example}\label{exY}
Let $m=5$ and $k=3$. Then
\[
Y_3=\left\{\begin{pmatrix}1&0&*&*&*\\0&1&*&*&*\\0&0&1&*&*\\0&0&0&1&*\\0&0&0&0&1\end{pmatrix}\right\},\qquad \psi_3(y)=\theta(y_{34}+y_{45}).
\]
Here $R_3=\{h_1h_2\}$, where $h_1=h_w\in\GL_4$ for a nonzero $w\in\FF_q^4$ and $h_2=h_{w'}\in\GL_3$ for a nonzero $w'\in\FF_q^3$, so $|R_3|=(q^4-1)(q^3-1)$. The actual elements of $R_3$ depend on the choices of the $h_w$.
\end{example}

\begin{example}\label{exW}
Let $\nu=(4,3,1)$ (so $m=8$) and $\mu=(3,1,1)$. The skew diagram $\nu/\mu$ (the shaded boxes below) consists of the boxes $(1,4)$, $(2,2)$, $(2,3)$, no two of which lie in the same column, so it is a horizontal $3$-strip and the strip letters are $6,7,8$. The row-reading tableaux of $\nu$ and $\mu$ are
\[
T_\nu=\ytab{1\&2\&3\&|[fill=gray!35]|4\\5\&|[fill=gray!35]|6\&|[fill=gray!35]|7\\8},\qquad T_\mu=\ytab{1\&2\&3\\4\\5}.
\]
The tableau $T'(\nu,\mu)$ agrees with $T_\mu$ on the boxes of $\mu$ and carries $6,7,8$ on the strip boxes, from left to right:
\[
T'(\nu,\mu)=\ytab{1\&2\&3\&|[fill=gray!35]|8\\4\&|[fill=gray!35]|6\&|[fill=gray!35]|7\\5}.
\]
Reading off the labels box by box, $W(\nu,\mu)$ sends $v_a$ to the vector whose index is the label of $T'(\nu,\mu)$ in the box where $T_\nu$ has $a$:
\[
W(\nu,\mu):\ v_4\mapsto v_8,\ v_5\mapsto v_4,\ v_8\mapsto v_5,\qquad v_1,v_2,v_3,v_6,v_7\ \text{fixed}.
\]
Thus the vectors in the boxes of the strip go to the strip letters $6,7,8$, and the vectors in the boxes of $\mu$ are renumbered $1,\dots,5$ by the row-reading order of $T_\mu$.
\end{example}

\section{The basis}\label{secbasis}

In Sections~\ref{secbasis}--\ref{secstd} the coefficient field is $\KK=\CC$, so $M^\lambda=M^\lambda_\CC$ and $S^\lambda=S^\lambda_\CC$; other fields are treated in Section~\ref{secmodular}.

\subsection{The theorem}

\begin{definition}\label{defB}
Put $B(\emptyset)=\{1\}$ and, for $\nu\vdash m\ge1$,
\[
B(\nu)=\bigl\{\,g\;\beta\;\Omega_{Y_k}\;W(\nu,\mu)\ :\ \nu/\mu\ \text{a horizontal $k$-strip},\ \beta\in B(\mu),\ g\in R_k\bigr\}\ \subset\ \CC[\GL_m].
\]
\end{definition}

\begin{theorem}[Main theorem]\label{thmmain}
For every partition $\lambda$ of $n$ and every prime power $q$, the vectors $\{\beta\,x_{T_\lambda}:\beta\in B(\lambda)\}$ form a basis of $S^\lambda$. In particular $|B(\lambda)|=\dim S^\lambda$.
\end{theorem}

Unrolling the recursion gives a closed form indexed by semistandard tableaux. Let $\tau$ be a semistandard tableau of shape $\lambda$ and content $\alpha=(\alpha_1,\dots,\alpha_r)$, with chain $\lambda=\nu^0\supset\nu^1\supset\dots\supset\nu^r=\emptyset$ (Section~\ref{secback}). Put $m_j=|\nu^j|$ and $k_j=m_j-m_{j+1}=\alpha_{r-j}$, and set
\[
W_j=W(\nu^j,\nu^{j+1}),\qquad E_j=\Omega_{Y^{(m_j)}_{k_j}},\qquad x_\tau=E_{r-1}W_{r-1}\cdots E_1W_1\,E_0W_0\;x_{T_\lambda}.
\]
(The factors do not commute in general: $E_0W_0$ acts first, and the later factors lie in the subgroup $\GL_{m_1}$ acting on the first $m_1$ coordinates.)

\begin{corollary}[Closed form]\label{corclosed}
For a composition $\alpha=(\alpha_1,\dots,\alpha_r)$ of $n$ put $s_i=\alpha_1+\dots+\alpha_i$ and
\[
w_\alpha(q)=\prod_{\substack{1\le i\le n-1\\ i\notin\{s_1,\dots,s_{r-1}\}}}(q^i-1).
\]
Then $S^\lambda$ has the basis
\[
\bigsqcup_{\alpha}\ \bigsqcup_{\tau}\ \{\,g_0g_1\cdots g_{r-1}\,x_\tau:\ g_j\in R^{(m_j)}_{k_j}\},
\]
where $\alpha$ runs over the compositions of $n$ and $\tau$ over the semistandard tableaux of shape $\lambda$ and content $\alpha$. Each $\tau$ contributes $w_\alpha(q)$ vectors, so
\[
\dim S^\lambda=\sum_{\alpha}K_{\lambda\alpha}\,w_\alpha(q),
\]
where $K_{\lambda\alpha}$ is the number of semistandard tableaux of shape $\lambda$ and content $\alpha$.
\end{corollary}

\begin{remark}
The vector $x_\tau$ is a $\psi$-eigenvector for each of the groups $Y^{(m_j)}_{k_j}$, since the later factors lie in $\GL_{m_{j+1}}$, which normalizes $Y^{(m_j)}_{k_j}$ and fixes $\psi_{k_j}$. The dimension formula of the corollary is the ``derivative'' recursion $\dim S^\lambda=\sum_\mu\dim S^\mu\prod_{i=n-k+1}^{n-1}(q^i-1)$ of Proposition~\ref{proplayer}. For the determinantal formula for $\dim S^\lambda$ (in James's labeling), this recursion is \cite[Thm.~10.21]{James84}.
\end{remark}

\begin{example}\label{exSt}
\textbf{The Steinberg module $S^{(2)}$} ($n=2$, $\dim=q$). For $\lambda=(2)$ the flags in $X_{(2)}$ are the flags $\FF_q^2\supsetneq W\supsetneq0$ with $W$ a line. There are $q+1$ of them: $F_a$ for $a\in\FF_q$, with $W=\mathrm{span}\{v_2+av_1\}$, and $F_\infty$, with $W=\mathrm{span}\{v_1\}$. Put
\[
T_a=\ytab{v_1\&\;v_2+av_1\;},\qquad T_\infty=\ytab{v_2\&v_1}.
\]
Then $F_{T_a}=F_a$ and $F_{T_\infty}=F_\infty$, so $m_{T_a}=[F_a]$ and $m_{T_\infty}=[F_\infty]$. In particular $T_\lambda=T_0$ and $x_{T_\lambda}=[F_0]-[F_\infty]=m_{T_0}-m_{T_\infty}$ (Example~\ref{exxT}). The group $Y_1=Y_2=\{y_a\}$, where $y_a:\ v_2\mapsto v_2+av_1$, satisfies $y_aF_0=F_a$ and fixes $F_\infty$. The two semistandard tableaux of shape $(2)$ give the following vectors.
\begin{center}
\renewcommand{\arraystretch}{1.5}\setlength{\tabcolsep}{9pt}
\begin{tabular}{l|c|c}
$\tau$ & $\ytab{1\&2}$ & $\ytab{1\&1}$\\\hline
content $\alpha$ & $(1,1)$ & $(2)$\\
$k$ & $1$ & $2$\\
$\psi_k(y_a)$ & $1$ & $\theta(a)$\\
$R_k$ & $\{1\}$ & $\{\mathrm{diag}(c,1):\ c\in\FF_q^\times\}$\\
basis vectors & $b_\infty=\Omega_{Y_1}x_{T_\lambda}$ & $b_c=\mathrm{diag}(c,1)\,\Omega_{Y_2}x_{T_\lambda}$\\
as flags & $b_\infty=\sum_{a}[F_a]-q\,[F_\infty]$ & $b_c=\sum_a\theta(-c^{-1}a)\,[F_a]$\\
as tableaux & $b_\infty=\sum_{a}m_{T_a}-q\,m_{T_\infty}$ & $b_c=\sum_a\theta(-c^{-1}a)\,m_{T_a}$\\
number of vectors & $1$ & $q-1$
\end{tabular}
\end{center}
\end{example}

\begin{example}\label{ex21}
\textbf{The module $S^{(2,1)}$} ($n=3$, $\dim=q^2+q$). For $\lambda=(2,1)$ the flags in $X_{(2,1)}$ are the flags $\FF_q^3\supsetneq L\supsetneq0$ with $L$ a line, so $M^{(2,1)}=\CC[\mathbb P^2(\FF_q)]$ and $S^{(2,1)}$ is the space of combinations whose coefficients sum to $0$. Let $\ell=\mathrm{span}\{v_1,v_2\}$. The flags with $L\subset\ell$ are the flags $F_a$ ($a\in\FF_q$), with $L=\mathrm{span}\{v_2+av_1\}$, and $F_\infty$, with $L=\mathrm{span}\{v_1\}$. The other flags are $F_{a,b}$ ($a,b\in\FF_q$), with $L=\mathrm{span}\{v_3+av_1+bv_2\}$. Put
\[
T_a=\ytab{v_1\&\;v_2+av_1\;\\v_3},\qquad T_\infty=\ytab{v_2\&v_1\\v_3},\qquad T_{a,b}=\ytab{v_1\&\;v_3+av_1+bv_2\;\\v_2}.
\]
Then $F_{T_a}=F_a$, $F_{T_\infty}=F_\infty$, and $F_{T_{a,b}}=F_{a,b}$, so $m_{T_a}=[F_a]$, $m_{T_\infty}=[F_\infty]$, and $m_{T_{a,b}}=[F_{a,b}]$. In particular $T_\lambda=T_0$, and since $N(T_\lambda)$ fixes $F_0$ we get $x_{T_\lambda}=q^2\bigl([F_0]-[F_\infty]\bigr)$.

Let $y_a:\ v_2\mapsto v_2+av_1$ and $y_{a,b}:\ v_3\mapsto v_3+av_1+bv_2$. Then $Y^{(3)}_1=\{y_{a,b}\}$ with $\psi_1=1$; $Y^{(2)}_1=Y^{(2)}_2=\{y_a\}$ with $\psi_1=1$ and $\psi_2(y_a)=\theta(a)$; and $Y^{(3)}_2$ is the group of upper unitriangular matrices with $\psi_2(y)=\theta(y_{23})$. The only permutation matrix other than $1$ that occurs is the matrix $s$ swapping $v_2$ and $v_3$, which is $W(\nu,\mu)$ for $\nu=(2,1)$ and $\mu=(1,1)$ or $\mu=(1)$. The semistandard tableaux of shape $(2,1)$ have content $(1,1,1)$ (two tableaux), $(2,1)$, or $(1,2)$, so $\dim S^{(2,1)}=2\cdot1+(q-1)+(q^2-1)$. These four tableaux give the following vectors.
\begin{center}
\renewcommand{\arraystretch}{1.5}\setlength{\tabcolsep}{9pt}
\begin{tabular}{l|c|c}
$\tau$ & $\ytab{1\&2\\3}$ & $\ytab{1\&1\\2}$\\\hline
content $\alpha$ & $(1,1,1)$ & $(2,1)$\\
$R$ & $\{1\}$ & $R^{(2)}_2=\{\mathrm{diag}(c,1):\ c\in\FF_q^\times\}$\\
basis vectors & $b_\infty=\Omega_{Y^{(2)}_1}\Omega_{Y^{(3)}_1}x_{T_\lambda}$ & $b_c=\mathrm{diag}(c,1)\,\Omega_{Y^{(2)}_2}\Omega_{Y^{(3)}_1}x_{T_\lambda}$\\
as flags & $b_\infty=q^4\bigl(\sum_{a}[F_a]-q\,[F_\infty]\bigr)$ & $b_c=q^4\sum_a\theta(-c^{-1}a)\,[F_a]$\\
as tableaux & $b_\infty=q^4\bigl(\sum_{a}m_{T_a}-q\,m_{T_\infty}\bigr)$ & $b_c=q^4\sum_a\theta(-c^{-1}a)\,m_{T_a}$\\
number of vectors & $1$ & $q-1$
\end{tabular}
\end{center}
\medskip
\begin{center}
\renewcommand{\arraystretch}{1.5}\setlength{\tabcolsep}{9pt}
\begin{tabular}{l|c|c}
$\tau$ & $\ytab{1\&3\\2}$ & $\ytab{1\&2\\2}$\\\hline
content $\alpha$ & $(1,1,1)$ & $(1,2)$\\
$R$ & $\{1\}$ & $R^{(3)}_2=\{h_w:\ w\in\FF_q^2\setminus0\}$\\
basis vectors & $b'_\infty=\Omega_{Y^{(2)}_1}\Omega_{Y^{(3)}_1}\,s\,x_{T_\lambda}$ & $b_w=h_w\,\Omega_{Y^{(3)}_2}\,s\,x_{T_\lambda}$\\
as flags & $b'_\infty=q^3\bigl(\sum_{a,b}[F_{a,b}]-q^2\,[F_\infty]\bigr)$ & $b_w=q^3\sum_{a,b}\theta(-w_1a-w_2b)\,[F_{a,b}]$\\
as tableaux & $b'_\infty=q^3\bigl(\sum_{a,b}m_{T_{a,b}}-q^2\,m_{T_\infty}\bigr)$ & $b_w=q^3\sum_{a,b}\theta(-w_1a-w_2b)\,m_{T_{a,b}}$\\
number of vectors & $1$ & $q^2-1$
\end{tabular}
\end{center}
\end{example}

\subsection{Outline of the proof}

The proof of Theorem~\ref{thmmain} is an induction on $n$ built on the following four steps.
\begin{enumerate}
\item (Section~\ref{seclayers}) Let $X$ be any module of the mirabolic subgroup $P_n$. Then, as a vector space, $X$ is the direct sum of the translates $g\,\Omega_{Y_k}X$ ($1\le k\le n$, $g\in R_k$) of its \emph{layers} $\Omega_{Y_k}X$, and each layer is a $\GL_{n-k}$-module (Lemma~\ref{lemlayers}).
\item (Section~\ref{seclayersS}) For $X=S^\lambda$ the layer $M_k=\Omega_{Y_k}S^\lambda$ is $\bigoplus_{\lambda/\mu}S^\mu$ as a $\GL_{n-k}$-module, the sum over horizontal $k$-strips (Proposition~\ref{proplayer}).
\item (Sections~\ref{secinner}--\ref{secdiag}) The vector $\Omega_{Y_k}x_{T'(\lambda,\mu)}$ has a nonzero component in the summand $M^\mu_k\cong S^\mu$ (Theorem~\ref{thmD1}), no component in $M^{\mu'}_k$ unless $\mu'\trianglerighteq\mu$ (Theorem~\ref{thmT}), and its $\mu$-component is the image of $x_{T_\mu}$ under a fixed isomorphism (Theorem~\ref{thmD2}).
\item (Section~\ref{secproof}) If $\{\beta x_{T_\mu}:\beta\in B(\mu)\}$ is a basis of $S^\mu$ for each $\mu$ with $\lambda/\mu$ a horizontal $k$-strip, then $\{\beta \Omega_{Y_k}x_{T'(\lambda,\mu)}\}$ is a basis of $M_k$.
\end{enumerate}
These steps prove the theorem by induction on $n$ (Section~\ref{secproof}).

\section{The vector \texorpdfstring{$x_T$}{xT}}\label{secx}

In this section we prove the properties of $x_T$ that are used later. It suffices to consider $T=T_\lambda$ (for $T=gT_\lambda$ everything is transported by $g$, by the definitions of Section~\ref{subflags}). Write $F=F_T$, $P=P(T)$, $N=N(T)$, and let
\[
Q=\{g:\ g_{ab}\ne0\Rightarrow\row(a)\le\row(b)\},\qquad L=\{g:\ g_{ab}\ne0\Rightarrow\row(a)=\row(b)\}\cong\prod_i\GL_{\lambda_i}(\FF_q).
\]
Then $Q=L\ltimes N$ is the parabolic subgroup of type $\lambda$ with Levi factor $L$ and unipotent radical $N$, $P$ is a parabolic subgroup of type $\lambda'$, $R(T)=W_L$ is the Weyl group of $L$, and $B_L:=L\cap P$ is a Borel subgroup of $L$ containing the diagonal torus, with $W_L\cap B_L=1$.

\begin{lemma}\label{lemfact}
If $l_1,l_2\in L$, $n_1,n_2\in N$, and $l_1n_1F=l_2n_2F$, then $l_1B_L=l_2B_L$.
\end{lemma}
\begin{proof}
The element $z=(l_2n_2)^{-1}l_1n_1$ fixes $F$, so $z\in P$, and $z\in Q$. Since $Q=L\ltimes N$, the image of $z$ in $Q/N\cong L$ is its block diagonal part, which is $l_2^{-1}l_1$ because $z\equiv l_2^{-1}l_1$ modulo $N$. The nonzero entries of the block diagonal part of $z$ are among the nonzero entries of $z$, so they lie in positions allowed for $P$. Hence $l_2^{-1}l_1\in P\cap L=B_L$.
\end{proof}

\begin{proposition}\label{propstpull}
Recall from Section~\ref{subflags} that $x(F')$ denotes the coefficient of $[F']$ in $x\in M^\lambda$. Then:
\begin{enumerate}
\item[(i)] The sets $NwF$ ($w\in W_L$) are pairwise disjoint, $x_T(F')=|N\cap P|\,\sgn(w)$ for $F'\in NwF$, and $x_T(F')=0$ for $F'\notin\bigcup_wNwF$. In particular $x_T\ne0$.
\item[(ii)] $N$ fixes $x_T$, and $\CC[L]x_T\cong\St_L$.
\end{enumerate}
\end{proposition}
\begin{proof}
Let $y=\sum_{F'\in NF}[F']$. The stabilizer of $F$ in $N$ is $N\cap P$, so $\sum_{n\in N}n\,[F]=|N\cap P|\,y$. Each $w\in R(T)=W_L$ normalizes $N$, so $w\,y=\sum_{F'\in NwF}[F']$. By the definition of $x_T$,
\[
x_T=\sum_{w\in W_L}\sgn(w)\,w\sum_{n\in N}n\,[F]=|N\cap P|\sum_{w\in W_L}\sgn(w)\,w\,y.
\]
If $w\,n_1F=w'n_2F$ with $w,w'\in W_L$ and $n_1,n_2\in N$, then $wB_L=w'B_L$ by Lemma~\ref{lemfact}, so $w=w'$ since $W_L\cap B_L=1$. Hence the sets $NwF$ are pairwise disjoint, and the formula gives the values of $x_T$ in (i); since $|N\cap P|\ne0$, $x_T\ne0$. Each $w\,y$ is a sum over an $N$-orbit, so $N$ fixes $w\,y$ and hence $x_T$.

Finally, let $L$ act on $\CC[L/B_L]$ by left multiplication. Since $B_L$ fixes $y$, the rule $lB_L\mapsto l\,y$ is a well-defined $L$-map $\CC[L/B_L]\to M^\lambda$. It is injective: by Lemma~\ref{lemfact} the vectors $l\,y$ for distinct cosets have pairwise disjoint supports $lNF$, so they are linearly independent. Hence $lB_L\mapsto l\,y$ is an isomorphism $\CC[L/B_L]\to\CC[L]y$. It sends $wB_L$ to $w\,y$ for $w\in W_L$, so by the formula for $x_T$ it sends $\sum_{w\in W_L}\sgn(w)\,wB_L$ to $|N\cap P|^{-1}x_T$. Thus $\CC[L]x_T$ is isomorphic to the submodule of $\CC[L/B_L]$ generated by $\sum_{w\in W_L}\sgn(w)\,wB_L$. By Steinberg's theorem~\cite{St}, this submodule is the (irreducible) Steinberg module of $L$. So $\CC[L]x_T\cong\St_L$.
\end{proof}

\begin{theorem}\label{thmx}
Let $T$ be any $\FF_q^n$-tableau of shape $\lambda$. Then $x_T\in S^\lambda$, $x_T\ne0$, and $\CC[G]\,x_T=S^\lambda$.
\end{theorem}

\begin{proof}
Let $V=\CC[G]x_T\subseteq M^\lambda$.

\emph{Step 1: every irreducible constituent of $V$ is $S^\nu$ with $\nu\trianglerighteq\lambda$.} By Proposition~\ref{propstpull}, $\CC[Q]x_T=\CC[L]x_T\cong\Inf_L^Q\St_L$ (the inflation of $\St_L$ to $Q=L\ltimes N$, with $N$ acting trivially). The inclusion $\CC[Q]x_T\subseteq M^\lambda$ is a $Q$-map, so by Frobenius reciprocity it gives a $G$-map $\Ind_Q^G\Inf_L^Q\St_L\to M^\lambda$, $g\otimes v\mapsto gv$, whose image is $\CC[G]x_T=V$. Thus $V$ is a quotient of $\Ind_Q^G\Inf_L^Q\St_L$, and by (H1) the character of this module is $\sum_{\nu\trianglerighteq\lambda}K_{\nu\lambda}\chi^\nu$.

\emph{Step 2: $V=S^\lambda$.} By Section~\ref{subflags}(e), every constituent of $M^\lambda$ is $S^\nu$ with $\nu\domeq\lambda$, and $S^\lambda$ occurs exactly once. Together with Step~1 this shows that every constituent of $V$ is $S^\lambda$. Since $V\ne0$ and $S^\lambda$ occurs only once in $M^\lambda$, $V=S^\lambda$.
\end{proof}

\section{Layers of modules of the mirabolic subgroup}\label{seclayers}

For $n\ge1$ let $N_n=Y^{(n)}_1=\{t_a=1+\sum_{i<n}a_iE_{in}\}\cong\FF_q^{n-1}$ be the group of matrices that differ from the identity only in the first $n-1$ entries of the last column, and let $P_n=\GL_{n-1}\ltimes N_n$ be the \emph{mirabolic subgroup} (the matrices of $\GL_n$ with last row $e_n^T$; see Example~\ref{exmirab}); here $ht_ah^{-1}=t_{ha}$ for $h\in\GL_{n-1}$. For a row vector $\ell\in\FF_q^{n-1}$ put
\[
e_\ell=|N_n|^{-1}\sum_{a}\theta(-\ell a)\,t_a\in\CC[N_n].
\]
These are orthogonal idempotents with sum $1$, and $t_ae_\ell=\theta(\ell a)e_\ell$, $h\,e_\ell\,h^{-1}=e_{\ell h^{-1}}$ for $h\in\GL_{n-1}$. Let $\ell_1=(0,\dots,0,1)$. Its stabilizer in $\GL_{n-1}$ is $P_{n-1}\supset U_{n-1}$, and $\GL_{n-1}$ is transitive on the nonzero $\ell$; the matrices $h_w\in\GL_{n-1}$ of Section~\ref{subnotation}, where $w\in\FF_q^{n-1}\setminus0$, satisfy $h_we_{\ell_1}h_w^{-1}=e_w$.

\begin{example}\label{exmirab}
The mirabolic subgroup $P_n$ is the set of invertible $n\times n$ matrices whose last row is $(0,\dots,0,1)$. For $n=4$ the matrices in $P_4$ look like
\[
\left(\begin{array}{ccc|c}
\ast&\ast&\ast&a_1\\
\ast&\ast&\ast&a_2\\
\ast&\ast&\ast&a_3\\\hline
0&0&0&1
\end{array}\right)
=\begin{pmatrix}A&a\\0&1\end{pmatrix},\qquad A\in\GL_3,\ a\in\FF_q^3,
\]
where the upper left block $A$ is any invertible matrix, the last column $a$ is arbitrary, and only the bottom row is constrained. So $|P_n|=q^{n-1}|\GL_{n-1}|$. The product is
\[
\begin{pmatrix}A&a\\0&1\end{pmatrix}\begin{pmatrix}B&b\\0&1\end{pmatrix}=\begin{pmatrix}AB&Ab+a\\0&1\end{pmatrix},
\]
so $P_n$ acts on the vectors $\binom{x}{1}$ by $x\mapsto Ax+a$: it is the group of affine transformations of $\FF_q^{n-1}$. Its normal subgroup $N_n$ is the set of matrices with $A=1$ (the translations), and $P_n/N_n\cong\GL_{n-1}$. For $n=2$,
$P_2=\Bigl\{\begin{pmatrix}c&b\\0&1\end{pmatrix}:\ c\in\FF_q^\times,\ b\in\FF_q\Bigr\}$ is the ``$ax+b$ group'' of order $q(q-1)$. The groups $Y^{(n)}_k$ of Section~\ref{subnotation} (upper unitriangular matrices whose off-diagonal entries lie in the last $k$ columns) lie in $P_n$.
\end{example}

Directly from the definitions, $Y^{(n)}_k=N_n\,Y^{(n-1)}_{k-1}$ (with unique factorization $t_ay'$), $\psi^{(n)}_k(t_ay')=\theta(a_{n-1})\psi^{(n-1)}_{k-1}(y')$ for $k\ge2$, and hence
\[
\Omega_{Y^{(n)}_1}=|N_n|\,e_0,\qquad \Omega_{Y_k^{(n)}}=|N_n|\,e_{\ell_1}\,\Omega_{Y_{k-1}^{(n-1)}}\quad(k\ge2).
\]
Also $R^{(n)}_k=H_nR^{(n-1)}_{k-1}$ for $k\ge2$, where $H_n=\{h_w\}\subset\GL_{n-1}$ is the set of the matrices $h_w$ of Section~\ref{subnotation} with $w\in\FF_q^{n-1}\setminus0$. The group $\GL_{n-k}$ normalizes $Y_k^{(n)}$ and fixes $\psi_k$, so it commutes with $\Omega_{Y^{(n)}_k}$.

\begin{lemma}[Layer decomposition]\label{lemlayers}
For every $\CC[P_n]$-module $X$,
\[
X=\bigoplus_{k=1}^n\ \bigoplus_{g\in R_k}g\cdot \Omega_{Y_k}X ,
\]
as vector spaces, and $\Omega_{Y_k}X$ is a $\GL_{n-k}$-module. Consequently $\dim X=\sum_{k=1}^n\prod_{j=n-k+1}^{n-1}(q^j-1)\cdot\dim \Omega_{Y_k}X$.
\end{lemma}
This is \cite[Thm.~9.11]{James84}, which is proved there by an induction of the same kind (for modules of the full group, with the idempotents $E_r$ corresponding to $\Omega_{Y_r}$).
\begin{proof}
We use induction on $n$; the case $n=1$ is trivial. Let $n\ge2$. The idempotents $e_\ell$ decompose $X=e_0X\oplus\bigoplus_{\ell\ne0}e_\ell X$. Put $Z=e_{\ell_1}X$. For $\ell\neq0$, $e_\ell X=h_\ell Z$, since $h_\ell e_{\ell_1}X=e_\ell h_\ell X=e_\ell X$. So
\[
X=e_0X\oplus\bigoplus_{h\in H_n}hZ,\qquad e_0X=\Omega_{Y_1}X.
\]
The space $Z$ is stable under $P_{n-1}$ (which fixes $\ell_1$), so by induction $Z=\bigoplus_{k'}\bigoplus_{g'\in R^{(n-1)}_{k'}}g'\Omega_{Y^{(n-1)}_{k'}}Z$, and $\Omega_{Y^{(n-1)}_{k'}}Z=\Omega_{Y^{(n-1)}_{k'}}e_{\ell_1}X=\Omega_{Y^{(n)}_{k'+1}}X$ up to the scalar $|N_n|$. Substituting gives the claim, because $R_k^{(n)}=H_nR_{k-1}^{(n-1)}$.
\end{proof}

\section{The layers of \texorpdfstring{$S^\lambda$}{S(lambda)}}\label{seclayersS}

Let $M_k=\Omega_{Y_k}S^\lambda$ with $Y_k=Y^{(n)}_k$; this is a $\GL_{n-k}$-module, the \emph{$k$-th layer} of $S^\lambda$.

\begin{proposition}\label{proplayer}
$M_k\cong\bigoplus_{\lambda/\mu\ \text{horizontal $k$-strip}}S^\mu$ as $\GL_{n-k}$-modules, and the sum is multiplicity free. We write $M_k^\mu$ for the summand isomorphic to $S^\mu$. Consequently
\[
\dim S^\lambda=\sum_{k=1}^n\ \prod_{j=n-k+1}^{n-1}(q^j-1)\ \sum_{\lambda/\mu}\dim S^\mu .
\]
\end{proposition}
This is the characteristic-$0$ case of~\cite[Thm.~16.9]{James84}, which gives a filtration of each layer by smaller Specht modules that splits in characteristic $0$; see also his Branching Theorem~\cite[Thm.~16.10]{James84}. James's proof uses his Kernel Intersection Theorem; the proof below uses Harish-Chandra theory instead.
\begin{proof}
The group $Y_k$ is the semidirect product $V\rtimes U_k$, where $V$ is the unipotent radical of the parabolic subgroup of type $(n-k,k)$, and $\psi_k$ is trivial on $V$. Hence $\Omega_{Y_k}=\Omega_{U_k,\psi}\,\Omega_V$ (a product of the unnormalized Fourier sums over $U_k$ and $V$, each a nonzero scalar multiple of the projection onto the corresponding eigenspace). Note that $\Omega_VS^\lambda$ is the Harish-Chandra restriction of $S^\lambda$ to $\GL_{n-k}\times\GL_k$ (the Jacquet module), and $\Omega_{U_k,\psi}$ projects it onto its $(U_k,\psi)$-eigenspace. So $M_k$ is the $(U_k,\psi)$-eigenspace of the Harish-Chandra restriction of $S^\lambda$. By Section~\ref{subunip}, (H1) and (H2), this is $\bigoplus_\mu c^\lambda_{\mu,(k)}S^\mu$. The dimension formula follows from Lemma~\ref{lemlayers}.
\end{proof}

By Lemma~\ref{lemlayers} and Proposition~\ref{proplayer}, a basis of $S^\lambda$ is obtained from bases of the layers $M_k$, and a basis of $M_k$ from bases of the $S^\mu$, once we have a way of embedding each $S^\mu$ in $M_k$ via a vector. The next four sections construct that vector; the idea of using such a decomposition goes back to Bernstein and Zelevinsky~\cite{BZ,BZ2} (for local fields); for finite general linear groups the corresponding operators, the decomposition of Lemma~\ref{lemlayers}, and the dimension formula appear in~\cite[Thm.~9.11, Cor.~10.16, Thm.~10.21]{James84} and in~\cite{Jam86}, Def.~1.7 and Lemma~1.8.

\section{Inner products with \texorpdfstring{$x_T$}{xT}}\label{secinner}

Fix a coordinate tableau $T$ of shape $\lambda$ and write $P=P(T)$, $N=N(T)$, $U=U(T)$, $F=F_T$, and $A=N\cap U$ (so $A$ consists of the matrices supported on the positions with $\row(a)<\row(b)$ and $\col(a)<\col(b)$). The lemma below says that, in the permutation module $M^\lambda$ on flags, the $N$-orbits of the flags $wF$ ($w\in R(T)$) meet the set $UF$ only in a controlled way.

\begin{lemma}\label{lemorbits}
We have the following:
\begin{enumerate}
\item[(a)] $N=A\,(N\cap P)$ and $A\cap(N\cap P)=1$; in particular $|N|=|A|\,|N\cap P|$.
\item[(b)] $\{u\in U:\ uF\in NF\}=A$.
\item[(c)] For $w\in R(T)$ with $w\ne1$, $N\,wF\cap UF=\emptyset$.
\end{enumerate}
\end{lemma}
\begin{proof}
(a) The groups $A$ and $N\cap P$ are subgroups of $N$ whose sets of positions (those of $N$ with $\col(a)<\col(b)$, respectively with $\col(a)\ge\col(b)$) partition the set of positions of $N$, so $|N|=|A|\,|N\cap P|$. The sets of positions defining $U$ and $P$ are disjoint, so $U\cap P=1$, and hence $A\cap(N\cap P)=1$. Therefore $N=A\,(N\cap P)$.

(b) Since $A=N\cap U$, we have $AF\subseteq NF$. Conversely, let $u\in U$ with $uF=nF$ and $n\in N$. Write $n=ap$ with $a\in A$ and $p\in N\cap P$ by (a). Then $nF=aF$, so $u^{-1}a\in U$ fixes $F$, i.e., $u^{-1}a\in U\cap P=1$, and $u=a\in A$.

(c) Let $\mathcal C_c=\{a:\col(a)\ge c\}$ and $\mathcal R_r=\{a:\row(a)=r\}$, and for a matrix $g$ let $g_c$ be its submatrix with rows and columns in $\mathcal C_c$.

\emph{Step 1: if $gF\in UF$, then $\det g_c\ne0$ for every $c$.} Write $g=up$ with $u\in U$ and $p\in P$. Let $b\in\mathcal C_c$. If $p_{xb}\ne0$, then $\col(x)\ge\col(b)\ge c$, so $x\in\mathcal C_c$. This has two consequences. First, for $a,b\in\mathcal C_c$ the sum $(up)_{ab}=\sum_xu_{ax}p_{xb}$ runs only over $x\in\mathcal C_c$, so $(up)_c=u_cp_c$. Second, $p_{xb}=0$ whenever $x\notin\mathcal C_c$ and $b\in\mathcal C_c$. So if the indices in $\mathcal C_c$ are listed before the others, $p$ has the block form
\[
p=\begin{pmatrix}p_c&*\\0&p_{c'}\end{pmatrix},
\]
where $p_{c'}$ is the submatrix of $p$ with rows and columns outside $\mathcal C_c$. Hence $\det p=\det p_c\cdot\det p_{c'}$, and $p_c$ is invertible because $\det p\ne0$. Also $u_c$ is unitriangular when $\mathcal C_c$ is ordered by column, since $u_{ab}\ne0$ with $a\ne b$ forces $\col(a)<\col(b)$. So $\det g_c=\det p_c\ne0$.

\emph{Step 2: if $g=nw$ with $n\in N$ and $w\in R(T)$, then $\det g_c\ne0$ only if $w(\mathcal C_c\cap\mathcal R_r)=\mathcal C_c\cap\mathcal R_r$ for every $r$.} We have $g_{ab}=n_{a,w(b)}$. Since $w$ preserves rows and $n_{aa'}\ne0$ with $a\ne a'$ forces $\row(a)<\row(a')$, the entry $g_{ab}$ can be nonzero only if $\row(a)\le\row(b)$. So $g_c$ is block upper triangular for the partition of $\mathcal C_c$ by rows of $T$. In the diagonal block for row $r$ we have $g_{ab}=\delta_{a,w(b)}$, which is a matrix with a $1$ at $(w(b),b)$ whenever $w(b)\in\mathcal C_c$ and zeros elsewhere. It is invertible if and only if $w(\mathcal C_c\cap\mathcal R_r)=\mathcal C_c\cap\mathcal R_r$.

Now suppose $nwF\in UF$. By Steps~1 and~2, $w(\mathcal C_c\cap\mathcal R_r)=\mathcal C_c\cap\mathcal R_r$ for all $c$ and $r$. Each box $a$ of $T$ is the unique element of $(\mathcal C_c\cap\mathcal R_r)\setminus(\mathcal C_{c+1}\cap\mathcal R_r)$ for $c=\col(a)$ and $r=\row(a)$, so $w$ fixes every $a$, and $w=1$.
\end{proof}

The next corollary computes the inner product of $x_T$ with a linear combination of flags $uF$, $u\in U$. We use it in the proof of Theorem~\ref{thmD1} to evaluate $\langle x_{T'},\eta_\mu\rangle$, and in Corollary~\ref{corkappa} to find the constant $\kappa$.

\begin{corollary}[Inner product formula]\label{corpair}
For a subset $\mathcal U\subseteq U$ and scalars $\phi(u)$,
\[
\Bigl\langle x_T,\ \sum_{u\in\mathcal U}\phi(u)[uF]\Bigr\rangle=|N\cap P|\sum_{u\in\mathcal U\cap A}\overline{\phi(u)}.
\]
\end{corollary}
\begin{proof}
Since $U\cap P=1$, the flags $uF$ ($u\in U$) are distinct, so the inner product is $\sum_{u\in\mathcal U}\overline{\phi(u)}\,x_T(uF)$. By Lemma~\ref{lemorbits}(c) the flag $uF$ lies in no orbit $NwF$ with $w\ne1$, so by Proposition~\ref{propstpull}(i), $x_T(uF)=|N\cap P|$ if $uF\in NF$ and $x_T(uF)=0$ otherwise. By Lemma~\ref{lemorbits}(b), $uF\in NF$ if and only if $u\in A$.
\end{proof}

\begin{corollary}[The constant $\kappa$]\label{corkappa}
For every $\FF_q^n$-tableau $T$ of shape $\lambda$ we have $k_Tx_T=\kappa\,e_T$ with $\kappa=|N(T)|=q^{\sum_{r<s}\lambda_r\lambda_s}$.
\end{corollary}
\begin{proof}
By (A2), $k_Tx_T=\kappa\,e_T$ for some $\kappa\in\CC$. We must show $\kappa=|N(T)|$. Since $\kappa$ does not depend on $T$, we only need to consider $T=T_\lambda$. As $k_T$ is self-adjoint and $e_T=k_T[F_T]$ has coefficient $1$ at $[F_T]$, we have $\kappa=\langle k_Tx_T,[F_T]\rangle=\langle x_T,e_T\rangle$. By the inner product formula, with $\phi=\psi_T^{-1}$ and $\psi_T$ trivial on $A$, this equals $|N\cap P|\,|A|=|N|=q^{\sum_{r<s}\lambda_r\lambda_s}$.
\end{proof}

\section{The canonical map and nonvanishing of the diagonal}\label{secD1}

Let $\lambda/\mu$ be a horizontal $k$-strip, $k=n-|\mu|$, and let $T'=T'(\lambda,\mu)$. Put $H=\mathrm{span}\{v_1,\dots,v_{n-k}\}$ and $E=\mathrm{span}\{v_{n-k+1},\dots,v_n\}$; the labels $n-k+1,\dots,n$ are the strip letters, and $E_c=\mathrm{span}\{v_b:\ b\ \text{a strip letter},\ \col(b)\ge c\}$. Here $M^\mu=\CC[X_\mu]$ has the flags of $H$ as a basis. For a flag $F=(W_c)$ of $H$ (with $W_1=H$ and $W_c=0$ for $c>\mu_1$) put $JF=(W_c\oplus E_c)$, a flag of $\FF_q^n$ (its first term is $H\oplus E=\FF_q^n$), and let $J\colon M^\mu\to M^\lambda$ be the linear map with $J[F]=[JF]$. The map $J$ is injective and $\GL_{n-k}$-equivariant, and $J[F_{T_\mu}]=[F_{T'}]$.

Define $\Phi_{\rm can}=\Omega_{Y_k}\,\pi_\lambda\,J|_{S^\mu}$, where $\pi_\lambda$ is the orthogonal projection of $M^\lambda$ onto $S^\lambda$ with respect to the Hermitian form of Section~\ref{subflags}(d). This is a $\GL_{n-k}$-map $S^\mu\to M_k$, because $\pi_\lambda$ is $G$-equivariant (the form is $G$-invariant and $S^\lambda$ is a submodule), $J$ is $\GL_{n-k}$-equivariant, and $\Omega_{Y_k}$ commutes with $\GL_{n-k}$. The vector $\Omega_{Y_k}x_{T'}\in M_k$ plays the role of the generator $x_{T_\mu}$ of $S^\mu$ in the layer $M_k$: its component in the summand $M_k^\mu\cong S^\mu$ is a nonzero multiple of $\Phi_{\rm can}(x_{T_\mu})$, and its other components lie in summands $M_k^\nu$ with $\nu\triangleright\mu$ (Theorems~\ref{thmD1}, \ref{thmT}, and~\ref{thmD2}).

\begin{theorem}[Nonvanishing of the diagonal]\label{thmD1}
$\Phi_{\rm can}$ is injective with image $M_k^\mu$, and
\[
\langle \Omega_{Y_k}x_{T'},\Phi_{\rm can}(e_{T_\mu})\rangle>0.
\]
In particular the $M_k^\mu$-component of $\Omega_{Y_k}x_{T'}$ is nonzero.
\end{theorem}
\begin{proof}
The map $\Phi_{\rm can}\colon S^\mu\to M_k$ is a homomorphism of $\GL_{n-k}$-modules, and $S^\mu$ is irreducible. So its image is $0$ or isomorphic to $S^\mu$. By Proposition~\ref{proplayer}, $M_k^\mu$ is the only submodule of $M_k$ isomorphic to $S^\mu$, so $\Phi_{\rm can}$ is either $0$ or an isomorphism onto $M_k^\mu$. Let $\eta_\mu=\Omega_{Y_k}J(e_{T_\mu})$. We have $\pi_\lambda \eta_\mu=\Phi_{\rm can}(e_{T_\mu})\in M_k^\mu$, and as $\Omega_{Y_k}x_{T'}\in S^\lambda$,
\[
\langle \Omega_{Y_k}x_{T'},\Phi_{\rm can}(e_{T_\mu})\rangle=\langle \Omega_{Y_k}x_{T'},\eta_\mu\rangle=|Y_k|\,\langle x_{T'},\eta_\mu\rangle,
\]
because $y\eta_\mu=\psi_k(y)\eta_\mu$ for $y\in Y_k$. We compute $\eta_\mu$. The group $U(T_\mu)\subset\GL_{n-k}$ normalizes $Y_k$ and fixes $\psi_k$, so
\[
\eta_\mu=\sum_{y\in Y_k,\ u\in U(T_\mu)}\psi_k(y)^{-1}\psi_{T_\mu}(u)^{-1}\,[u\,yF_{T'}].
\]
Let $Y^+=Y_k\cap U(T')$ and $Y^0=Y_k\cap P(T')$. Then $Y_k=Y^+Y^0$, because the positions of $Y_k$ are the disjoint union of those of $Y^+$ and $Y^0$ and $Y^+\cap Y^0\subseteq U(T')\cap P(T')=1$. The group $Y^0$ fixes $F_{T'}$, and $\psi_k|_{Y^0}=1$ (the character $\psi_k$ only involves entries $(i-1,i)$ with $i-1$ and $i$ strip letters, and these have $\col(i-1)<\col(i)$, so they vanish on $Y^0$). Write $y=y^+y^0$ with $y^+\in Y^+$ and $y^0\in Y^0$. Since $y^0$ fixes $F_{T'}$ and $\psi_k(y^0)=1$, each summand depends only on $y^+$, which gives the factor $|Y^0|$. For $u\in U(T_\mu)$ and $y\in Y^+$, the conjugate $uyu^{-1}$ has the same block of strip letters as $y$, and its entry $(a,b)$ with $a\in H$ is a combination of the entries $y_{a'b}$ with $\col(a)\le\col(a')<\col(b)$. So $U(T_\mu)\subset U(T')$ normalizes $Y^+$ and fixes $\psi_k$. It meets $Y^+$ trivially, so $\mathcal U^+=U(T_\mu)Y^+$ is a subgroup of $U(T')$, and $\psi^+(uy^+)=\psi_{T_\mu}(u)\,\psi_k(y^+)$ is a character of $\mathcal U^+$. Hence
\[
\eta_\mu=|Y^0|\sum_{\bar u\in\mathcal U^+}\psi^+(\bar u)^{-1}\,[\bar uF_{T'}].
\]
By the inner product formula (Corollary~\ref{corpair}, with $T=T'$),
\[
\langle x_{T'},\eta_\mu\rangle=|N'\cap P'|\,|Y^0|\sum_{\bar u\in\mathcal U^+\cap A'}\overline{\psi^+(\bar u)^{-1}},
\]
with $N'=N(T')$, $P'=P(T')$, $A'=N'\cap U(T')$. First, $A'\subseteq\mathcal U^+$. Let $g\in A'$. Each nonzero off-diagonal entry $(a,b)$ of $g$ has $\row(a)<\row(b)$ and $\col(a)<\col(b)$, so $a$ is not a strip letter. Hence $g=\begin{pmatrix}g_H&B\\0&1\end{pmatrix}$ with respect to $H\oplus E$, and $g=uy$ with $u=\mathrm{diag}(g_H,1)\in U(T_\mu)$ and $y=u^{-1}g\in Y_k\cap U(T')=Y^+$. Moreover $\psi^+(g)=\psi_{T_\mu}(u)\,\psi_k(y)=1$, since $g_{ab}=0$ for $a\ne b$ in a common row, and $y_{ab}=0$ for $a\ne b$ strip letters. So $\psi^+\equiv1$ on $A'=\mathcal U^+\cap A'$, the sum above equals $|A'|$, and since $|N'\cap P'|\,|A'|=|N'|$ (Lemma~\ref{lemorbits}(a)),
\[
\langle x_{T'},\eta_\mu\rangle=|Y^0|\,|N'|,\qquad\langle \Omega_{Y_k}x_{T'},\Phi_{\rm can}(e_{T_\mu})\rangle=|Y_k|\,|Y^0|\,|N'|>0.
\]
Hence $\Phi_{\rm can}\ne0$, so $\Phi_{\rm can}$ is injective with image $M_k^\mu$, and the $M_k^\mu$-component of $\Omega_{Y_k}x_{T'}$ is nonzero.
\end{proof}

\section{Triangularity}\label{sectri}

The vector $\Omega_{Y_k}x_{T'(\mu)}$, where $T'(\mu)=T'(\lambda,\mu)$, has nonzero component in $M_k^\mu$, but it also has components in other summands $M_k^{\nu}$ of the same layer. The next theorem shows that these other components are in summands that are higher in the dominance order, so that the matrix of the vectors $\beta \Omega_{Y_k}x_{T'(\mu)}$ with respect to the decomposition $M_k=\bigoplus_\mu M_k^\mu$ is block triangular.

Number the boxes of the strip $\lambda/\mu$ from left to right, $j=1,\dots,k$: the box $j$ has column $c_j(\mu)$ (strictly increasing in $j$) and row $r_j(\mu)$ (weakly decreasing), and carries the letter $n-k+j$ in $T'(\mu)$. Let $s_c(\mu)$ be the number of boxes of the strip in columns $\ge c$, and $t_r(\mu)=\sum_{i<r}(\lambda_i-\mu_i)$ the number of boxes of the strip in rows $<r$.

\begin{theorem}[Triangularity]\label{thmT}
If the $M_k^{\nu}$-component of $\Omega_{Y_k}x_{T'(\mu)}$ is nonzero, then $\nu\trianglerighteq\mu$. Hence $\beta \Omega_{Y_k}x_{T'(\mu)}\in\bigoplus_{\nu\trianglerighteq\mu}M_k^{\nu}$ for every $\beta\in\CC[\GL_{n-k}]$.
\end{theorem}
\begin{proof}
\emph{Reduction to supports.} By Theorem~\ref{thmD1} for $\nu$, $M_k^{\nu}$ is spanned by the vectors $g\,\pi_\lambda \eta_\nu$ ($g\in\GL_{n-k}$), where $\eta_\nu=|Y^0|\sum_{\bar u\in\mathcal U^+(\nu)}\psi^+(\bar u)^{-1}[\bar uF_{T'(\nu)}]$. Since $\Omega_{Y_k}x_{T'(\mu)}\in S^\lambda$, we have $\langle \Omega_{Y_k}x_{T'(\mu)},g\pi_\lambda \eta_\nu\rangle=|Y_k|\,\langle x_{T'(\mu)},g\eta_\nu\rangle$. So it suffices to show that if $\nu\not\trianglerighteq\mu$, then for every $g\in\GL_{n-k}$ the support $\bigcup_wN'wF'$ of $x_{T'(\mu)}$ (where $N'=N(T'(\mu))$, $F'=F_{T'(\mu)}$) does not meet $g\,U(T'(\nu))F_{T'(\nu)}$.

\emph{A rank condition.} Let $F=(W_c)$ be a flag in $g\,U(T'(\nu))F_{T'(\nu)}$, and fix $c$. Let $\Lambda^H_c=\mathrm{span}\{v_a:\ a\le n-k,\ \col_{T'(\nu)}(a)<c\}$ and $\Lambda^E_c=\mathrm{span}\{v_a:\ a>n-k,\ \col_{T'(\nu)}(a)<c\}$, and let $\mathfrak S_c=\{a>n-k:\ \col_{T'(\nu)}(a)\ge c\}$ be the set of strip letters of $T'(\nu)$ in columns $\ge c$, so $|\mathfrak S_c|=s_c(\nu)$. Then $W_c$ is complementary to $g\Lambda^H_c\oplus\Lambda^E_c$. Put $X=W_c\cap(H\oplus\Lambda^E_c)$. Both $X$ and $g\Lambda^H_c\oplus\Lambda^E_c$ lie in $H\oplus\Lambda^E_c$ (as $g\in\GL_{n-k}$), and they meet trivially. So $\dim X\le\dim(H\oplus\Lambda^E_c)-\dim(g\Lambda^H_c\oplus\Lambda^E_c)=(n-k)-\dim\Lambda^H_c$, the number of boxes of $\nu$ in columns $\ge c$. Let $P_c\colon W_c\to\mathrm{span}\{v_a:a\in\mathfrak S_c\}$ be the coordinate projection. Its kernel is $X$, since $H\oplus\Lambda^E_c=\mathrm{span}\{v_a:\ a\notin\mathfrak S_c\}$. Now $\dim W_c$ is the number of boxes of $\lambda$ in columns $\ge c$, and $\dim X$ is at most the number of boxes of $\nu$ in columns $\ge c$. So $\rank P_c=\dim W_c-\dim X$ is at least the number of boxes of $\lambda/\nu$ in columns $\ge c$, which is $s_c(\nu)$. Thus $P_c$ is surjective, since its target has dimension $|\mathfrak S_c|=s_c(\nu)$.

\emph{A support bound.} Now suppose that $F$ also has the form $nwF'$ with $n\in N'$ and $w\in R(T'(\mu))$, so $W_c=n\,\mathrm{span}\{v_{a'}:a'\in w\mathcal C_c\}$, where $\mathcal C_c$ is the set of letters in columns $\ge c$ of $T'(\mu)$. In the basis $nv_{a'}$ ($a'\in w\mathcal C_c$) of $W_c$ the matrix of $P_c$ (with respect to this basis of $W_c$ and the basis $v_a$, $a\in\mathfrak S_c$) has rows indexed by $a\in\mathfrak S_c$, columns indexed by $a'\in w\mathcal C_c$, and entries $n_{aa'}$ (the coefficient of $v_a$ in $nv_{a'}$). Here $n_{aa}=1$, and $n_{aa'}\ne0$ with $a\ne a'$ forces $\row(a)<\row(a')$ (rows of $T'(\mu)$).

Fix $r\ge1$ and let $\mathfrak S_{c,r}=\{a\in\mathfrak S_c:\ \row_{T'(\mu)}(a)\ge r\}$, the strip letters that lie in columns $\ge c$ of $T'(\nu)$ and in rows $\ge r$ of $T'(\mu)$. The rows of $P_c$ indexed by $\mathfrak S_{c,r}$ are supported on the columns $a'$ with $\row(a')\ge r$, and there are $\lambda_{\ge r,\ge c}:=\sum_{i\ge r}(\lambda_i-c+1)^+$ of these (row $i$ has $(\lambda_i-c+1)^+$ boxes in columns $\ge c$), since $w$ preserves rows. These rows are linearly independent, because $P_c$ is surjective, so
\[
|\mathfrak S_{c,r}|\le\lambda_{\ge r,\ge c}.
\]
To compute $|\mathfrak S_{c,r}|$, note that the letter $n-k+j$ lies in column $c_j(\nu)$ of $T'(\nu)$ and in row $r_j(\mu)$ of $T'(\mu)$, so $|\mathfrak S_{c,r}|=\#\{j:\ c_j(\nu)\ge c,\ r_j(\mu)\ge r\}$. Since $c_j(\nu)$ increases and $r_j(\mu)$ decreases with $j$, the first condition holds for the last $s_c(\nu)$ values of $j$ and the second for the first $k-t_r(\mu)$ values. Hence $|\mathfrak S_{c,r}|=(s_c(\nu)-t_r(\mu))^+$, and
\[
(s_c(\nu)-t_r(\mu))^+\le\lambda_{\ge r,\ge c}\qquad\text{for all $r$ and $c$.}
\]

\emph{Dominance.} Fix $r\ge1$ and put $c=\nu_r+1$. Since $\lambda/\nu$ is a horizontal strip, $\lambda_{r+1}\le\nu_r<c$, so $\lambda_{\ge r,\ge c}=\lambda_r-\nu_r$. For $i<r$ we have $\nu_i\ge\nu_r$, so the strip boxes in columns $\ge c$ are exactly those in rows $\le r$, and $s_c(\nu)=\sum_{i\le r}(\lambda_i-\nu_i)$. Hence $s_c(\nu)-t_r(\mu)=(\lambda_r-\nu_r)+\sum_{i<r}(\mu_i-\nu_i)$, and the bound $(s_c(\nu)-t_r(\mu))^+\le\lambda_r-\nu_r$ forces $\sum_{i<r}\mu_i\le\sum_{i<r}\nu_i$. This is true for all $r$, so $\nu\trianglerighteq\mu$.

The last sentence of the theorem holds because $\beta\in\CC[\GL_{n-k}]$ preserves each summand $M_k^{\nu}$.
\end{proof}

\section{The diagonal component is the image of \texorpdfstring{$x_{T_\mu}$}{x(Tmu)}}\label{secdiag}

Finally we identify the diagonal component of $\Omega_{Y_k}x_{T'}$, not only as nonzero, but as the image under $\Phi_{\rm can}$ of the generator $x_{T_\mu}$ of $S^\mu$. This is what allows us to apply the induction hypothesis to $\{\beta\,x_{T_\mu}\}$. Let $J^\dagger\colon M^\lambda\to M^\mu$ be the adjoint of $J$ for the Hermitian forms, and put $\xi=J^\dagger \Omega_{Y_k}\,x_{T'}\in M^\mu$. We have that
\[
\xi(F^H)=\sum_{y\in Y_k}\psi_k(y)\,x_{T'}(y\,JF^H)
\]
for every flag $F^H$ of $H$, by $\langle J^\dagger v,[F^H]\rangle=\langle v,[JF^H]\rangle$ and the substitution $y\mapsto y^{-1}$.

\begin{theorem}[Identification of the diagonal component]\label{thmD2}
Let $Y^0=Y_k\cap P(T')$ (see Theorem~\ref{thmD1}),
\[
D=\sum_{r<r'}\bigl(\lambda_r\lambda_{r'}-\mu_r\mu_{r'}\bigr),\qquad c=q^{D}\,|Y^0|,
\]
and let $\omega>0$ be defined by $\Phi_{\rm can}^\dagger\Phi_{\rm can}=\omega\cdot\mathrm{id}$. Then:
\begin{itemize}
\item[(i)] $\xi=c\,x_{T_\mu}$.
\item[(ii)] $c$ is a power of $q$, hence a unit in $A=\ZZ[\zeta_p,1/p]$.
\item[(iii)] $(\Omega_{Y_k}x_{T'})^\mu=c'\,\Phi_{\rm can}(x_{T_\mu})$ with $c'=|Y_k|\,c/\omega\ne0$.
\end{itemize}
\end{theorem}
\begin{proof}
\emph{Notation.} Write $T'=T'(\lambda,\mu)$ (for the remainder of the proof $\row$ and $\col$ will refer to $T'$), $N'=N(T')$, $P'=P(T')$, $F'=F_{T'}$, $\mathcal C_c=\{a:\col(a)\ge c\}$, and let $\mathfrak S_c$ be the set of strip letters in $\mathcal C_c$, so $|\mathfrak S_c|=s_c(\mu)$. Let $V\subset Y_k$ be the abelian group of matrices supported on the positions $(a,b)$ with $a\in H$ and $b$ a strip letter, and $V_{\rm up}=V\cap N'$ (the positions with $\row(a)<\row(b)$). Let $K$ and $N_E$ be the groups of matrices in $N'$ supported on the positions $(a,b)$ with $a\in H$ and with $a$ a strip letter, respectively. Then $K=V_{\rm up}N_\mu$, where $N_\mu=N(T_\mu)$ acts on $H$, and $N'=KN_E$ with $|N'|=|K||N_E|$. These are product decompositions of sets that come from disjoint sets of positions, and we use only the set equalities and the cardinalities. We will use the following fact about the strip in Steps 1, 2, and 4: for a strip box $(r,c)$ we have $c>\mu_r\ge\lambda_{r+1}$, so no box of $\lambda$ below row $r$ lies in a column $\ge c$.

Let $F^H=(W^H_c)$ be a flag of $H$ and $y\in Y_k$, and put $F=yJF^H=(W_c)$. Since $y$ is the identity on $H$ and unitriangular on $E$, $W_c=W^H_c\oplus yE_c$. Hence:
\begin{itemize}
\item[(F1)] the coordinate projection of $W_c$ onto $\mathrm{span}\{v_a:a\in\mathfrak S_c\}$ is onto, and
\item[(F2)] $W_c\cap H=W^H_c$ for all $c$, so $F^H$ is recovered from $F$ as $F\cap H$.
\end{itemize}

The vector $x_{T'}$ is supported on the sets $N'wF'$ with $w\in R(T')$. Step~1 uses (F1) to show that a flag $yJF^H$ can lie in $N'wF'$ only if $w\in R(T_\mu)$. Step~2 describes $N'wF'$ for those $w$. Step~3 uses (F2) to recover $F^H$, which identifies the support of $\xi$, and then computes $\xi$. Step~4 computes the constant.

\emph{Step 1: only $w\in R(T_\mu)$ occur.} Suppose a flag $F=yJF^H$ lies in $N'wF'$, with $w\in R(T')$. So $F=n'wF'$ with $n'\in N'$, and $F$ satisfies (F1). We show that $w\in R(T_\mu)$. In the basis $n'v_{a'}$, $a'\in w\mathcal C_c$, the matrix of the projection in (F1) has entries $n'_{aa'}$ ($a\in\mathfrak S_c$), where $n'_{aa}=1$ and $n'_{aa'}\ne0$ with $a\ne a'$ forces $\row(a)<\row(a')$. Let $b$ be a strip letter, $c=\col(b)$, $r=\row(b)$. By the fact about the strip, $w\mathcal C_c$ has no letters below row $r$, so if $b\notin w\mathcal C_c$ the row of $b$ in this matrix is zero, contradicting (F1). Hence $b\in w\mathcal C_{\col(b)}$ for every strip letter $b$, i.e., $w^{-1}(b)$ is a letter of row $r$ in a column $\ge\col(b)$, and so is a strip letter of row $r$ at or to the right of $b$. Thus $w^{-1}$ permutes the strip letters of each row and moves each of them weakly to the right, so it fixes them all. That is, $w\in R(T_\mu)$. Thus the only sets $N'wF'$ that meet the flags $yJF^H$ are those with $w\in R(T_\mu)$.

\emph{Step 2: for $w\in R(T_\mu)$, $N'wF'=V_{\rm up}N_\mu wF'$.} We show that $N_E$ fixes $wF'$. Let $a\ne b$ be a position with $n_{ab}\ne0$ for some $n\in N_E$. Then $a$ is a strip letter, with $r=\row(a)$, and $\row(b)>r$. As $w$ preserves rows and fixes the strip letters, in $wT'$ the letter $a$ is still in column $\col(a)$, while $b$ is in a column $\le\lambda_{\row(b)}\le\lambda_{r+1}<\col(a)$. So $\col_{wT'}(a)\ge\col_{wT'}(b)$ for every such position, i.e., $N_E\subseteq P(wT')=\Stab(wF')$. Hence $N'wF'=KN_EwF'=KwF'=V_{\rm up}N_\mu wF'$.

\emph{Step 3: the sum.} By definition $\xi(F^H)=\sum_{y\in Y_k}\psi_k(y)\,x_{T'}(yJF^H)$. If $yJF^H$ lies in the support of $x_{T'}$, say in $N'wF'$ with $w\in R(T')$, then Step 1 gives $w\in R(T_\mu)$ and Step 2 gives $yJF^H=v\,n_H\,wF'=vJ(n_HwF_\mu)$ with $v\in V_{\rm up}$, $n_H\in N_\mu$ (note $wF'=J(wF_\mu)$). Intersecting with $H$ and using (F2) on both sides gives $F^H=n_HwF_\mu$. So $\xi$ is supported on $\bigcup_{w\in R(T_\mu)}N_\mu wF_\mu$, the support of $x_{T_\mu}$.

Fix $F^H=gF_\mu$ in this set, with $g=n_Hw\in\GL_{n-k}$, so that $JF^H=gF'$. By the same argument, the $y\in Y_k$ with $yJF^H$ in the support of $x_{T'}$ are exactly those with $yJF^H\in V_{\rm up}JF^H$, i.e., the elements of $V_{\rm up}\Stab_{Y_k}(JF^H)$. For these $y$, $yJF^H\in N'wF'$, so $x_{T'}(yJF^H)=|N'\cap P'|\,\sgn(w)$ by Proposition~\ref{propstpull}(i). The element $g$ normalizes $Y_k$, and it normalizes $V_{\rm up}$ because it preserves each space $H_{<r}=\mathrm{span}\{v_a:a\in H,\ \row(a)<r\}$ ($n_H$ maps each row space into the span of itself and the earlier rows, and $w$ preserves each row). Hence $\Stab_{Y_k}(JF^H)=g(Y_k\cap P')g^{-1}=gY^0g^{-1}$ and $\Stab_{V_{\rm up}}(JF^H)=g(V_{\rm up}\cap P')g^{-1}$. The character $\psi_k$ is trivial on $V_{\rm up}$ (it is trivial on $V$) and on $gY^0g^{-1}$ (since $g$ fixes $\psi_k$ and $\psi_k|_{Y^0}=1$). Therefore $\xi(F^H)$ is $|N'\cap P'|\,\sgn(w)$ times the number of elements of $V_{\rm up}\,gY^0g^{-1}$, which is $|V_{\rm up}|\,|Y^0|/|V_{\rm up}\cap P'|$. That is,
\[
\xi(F^H)=|N'\cap P'|\,\sgn(w)\ \frac{|V_{\rm up}|\,|Y^0|}{|V_{\rm up}\cap P'|},
\]
and $\xi(F^H)=0$ for $F^H$ outside the support of $x_{T_\mu}$.

\emph{Step 4: the constant.} Comparing with $x_{T_\mu}(F^H)=|N_\mu\cap P_\mu|\sgn(w)$ (Proposition~\ref{propstpull}(i) for $T_\mu$), the sign cancels and
\[
c=\frac{\xi(F^H)}{x_{T_\mu}(F^H)}=\frac{|N'\cap P'|\,|V_{\rm up}|\,|Y^0|}{|V_{\rm up}\cap P'|\,|N_\mu\cap P_\mu|}.
\]
We evaluate this by classifying the positions $(a,b)$ of $N'\cap P'$ (those with $\row(a)<\row(b)$ and $\col(a)\ge\col(b)$) as follows. If $a,b\in H$, these are the positions of $N_\mu\cap P_\mu$. If $a$ is a strip letter, then every position of $N'$ is in $P'$ (this is the inequality $\col(a)>\lambda_{r+1}\ge\col(b)$ of Step 2 for $w=1$), so these are all the positions of $N_E$. If $a\in H$ and $b$ is a strip letter, these are the positions of $V_{\rm up}\cap P'$. Hence $|N'\cap P'|=|N_\mu\cap P_\mu|\,|N_E|\,|V_{\rm up}\cap P'|$, and
\[
c=|N_E|\,|V_{\rm up}|\,|Y^0|=\frac{|N'|}{|N_\mu|}\,|Y^0|=q^{D}\,|Y^0|,
\]
because $|N'|=q^{\sum_{r<r'}\lambda_r\lambda_{r'}}$, $|N_\mu|=q^{\sum_{r<r'}\mu_r\mu_{r'}}$, $|K||N_E|=|N'|$, and $|K|=|V_{\rm up}||N_\mu|$.

\emph{Step 5: proof of (iii).} The map $\Phi_{\rm can}\colon S^\mu\to M_k^\mu$ is an isomorphism of irreducible $\GL_{n-k}$-modules, so $\Phi_{\rm can}^\dagger\Phi_{\rm can}=\omega\,\mathrm{id}$ for some $\omega>0$ (Schur's lemma). Write $(\Omega_{Y_k}x_{T'})^\mu=\Phi_{\rm can}(s_0)$ with $s_0\in S^\mu$. For $s\in S^\mu$, using $\Omega_{Y_k}x_{T'}\in S^\lambda$ and the facts that $\Omega_{Y_k}$ is self-adjoint and satisfies $\Omega_{Y_k}^2=|Y_k|\,\Omega_{Y_k}$ (Section~\ref{subnotation}),
\[
\omega\langle s_0,s\rangle=\langle\Phi_{\rm can}s_0,\Phi_{\rm can}s\rangle=\langle \Omega_{Y_k}x_{T'},\Phi_{\rm can}s\rangle=\langle \Omega_{Y_k}x_{T'},\Omega_{Y_k}Js\rangle=|Y_k|\langle\xi,s\rangle=|Y_k|\,c\,\langle x_{T_\mu},s\rangle .
\]
Since $s\in S^\mu$ is arbitrary and $s_0,x_{T_\mu}\in S^\mu$, we get $s_0=(|Y_k|c/\omega)\,x_{T_\mu}$, i.e., $(\Omega_{Y_k}x_{T'})^\mu=(|Y_k|c/\omega)\,\Phi_{\rm can}(x_{T_\mu})$.
\end{proof}

\section{The lifting theorem and the proof of the main theorem}\label{secproof}

\begin{theorem}[Lifting]\label{thmlift}
For each horizontal $k$-strip $\lambda/\mu$ let $B_\mu\subset\CC[\GL_{n-k}]$ be such that $\{\beta x_{T_\mu}:\beta\in B_\mu\}$ is a basis of $S^\mu$. Then
\[
\{\,g\,\beta\,\Omega_{Y_k}\,x_{T'(\lambda,\mu)}:\ k\ge1,\ \lambda/\mu\ \text{a $k$-strip},\ \beta\in B_\mu,\ g\in R_k\}
\]
is a basis of $S^\lambda$.
\end{theorem}
\begin{proof}
By Lemma~\ref{lemlayers} applied to $X=S^\lambda|_{P_n}$, it suffices to show, for each $k$, that $\{\beta \Omega_{Y_k}x_{T'(\lambda,\mu)}\}$ is a basis of $M_k=\bigoplus_\mu M_k^\mu$. By Theorem~\ref{thmT}, $\beta \Omega_{Y_k}x_{T'(\lambda,\mu)}\in\bigoplus_{\nu\trianglerighteq\mu}M_k^{\nu}$. By Theorem~\ref{thmD2} its $M_k^\mu$-component is $\beta\,(\Omega_{Y_k}x_{T'(\lambda,\mu)})^\mu=c'\beta\,\Phi_{\rm can}(x_{T_\mu})=c'\,\Phi_{\rm can}(\beta x_{T_\mu})$, and since $\Phi_{\rm can}:S^\mu\to M_k^\mu$ is an isomorphism (Theorem~\ref{thmD1}) these components form a basis of $M_k^\mu$. Order the strips by a linear extension of the dominance order with the largest first. Then the matrix of the family with respect to $M_k=\bigoplus_\mu M_k^\mu$ is block triangular with invertible diagonal blocks, hence invertible.
\end{proof}

\begin{proof}[Proof of Theorem~\ref{thmmain}]
We use induction on $n$; $n=0$ is trivial ($S^\emptyset$ is one dimensional and $B(\emptyset)=\{1\}$). For $\lambda\vdash n$ and a $k$-strip $\mu$, the set $\{\beta x_{T_\mu}:\beta\in B(\mu)\}$ is a basis of $S^\mu$ by induction ($\mu\vdash n-k$ and $B(\mu)\subset\CC[\GL_{n-k}]$). By Theorem~\ref{thmlift} the vectors $g\beta \Omega_{Y_k}x_{T'(\lambda,\mu)}=g\beta \Omega_{Y_k}W(\lambda,\mu)x_{T_\lambda}$ form a basis of $S^\lambda$, and these are exactly the vectors $\beta'x_{T_\lambda}$ with $\beta'\in B(\lambda)$.
\end{proof}

\begin{proof}[Proof of Corollary~\ref{corclosed}]
Unrolling Definition~\ref{defB}: $\beta=g_0\beta_1E_0W_0$ with $\beta_1=g_1\beta_2E_1W_1$, and so on, which gives $\beta=g_0g_1\cdots g_{r-1}\,E_{r-1}W_{r-1}\cdots E_0W_0$. Finally, $|R^{(m_j)}_{k_j}|=\prod_{i=m_{j+1}+1}^{m_j-1}(q^i-1)$, and the integers $m_{j+1}<i<m_j$ ($0\le j<r$) are exactly the integers $1\le i\le n-1$ other than $m_1,\dots,m_{r-1}=s_{r-1},\dots,s_1$; so the product over $j$ is $w_\alpha(q)$.
\end{proof}

\section{Grouping by standard tableaux}\label{secstd}

Let $\tau$ be a semistandard tableau of shape $\lambda$ and content $\alpha$, and let $s=\std(\tau)$ be its standardization (Section~\ref{secback}). For a standard tableau $s$ put
\[
A(s)=\{\,i:\ i+1\ \text{lies in a column strictly to the right of}\ i\,\},\qquad \maj(s^T)=\sum_{i\in A(s)}i,
\]
the major index of the transposed tableau. Recall that $\tau$ contributes $w_\tau(q)=w_\alpha(q)$ vectors to the basis of Corollary~\ref{corclosed}.

\begin{proposition}[Grouping by standard tableaux]\label{propstd}
For every standard tableau $s$ of shape $\lambda$,
\[
\sum_{\tau:\ \std(\tau)=s}w_\tau(q)\;=\;q^{\maj(s^T)} .
\]
Hence the basis of Corollary~\ref{corclosed} is the disjoint union, over the standard tableaux $s$ of shape $\lambda$, of sets $\mathcal B_s$ with $|\mathcal B_s|=q^{\maj(s^T)}$, where $\mathcal B_s$ consists of the vectors $g_0\cdots g_{r-1}x_\tau$ with $\std(\tau)=s$; and $\dim S^\lambda=\sum_sq^{\maj(s^T)}$.
\end{proposition}
\begin{proof}
Fix $s$. For a semistandard tableau $\tau$ with $\std(\tau)=s$ put
\[
S(\tau)=\{i:\ \text{the labels } i,i+1 \text{ of } s \text{ receive the same entry of } \tau\}.
\]
We show that $\tau\mapsto S(\tau)$ is a bijection from $\{\tau:\ \std(\tau)=s\}$ onto the subsets of $A(s)$. (This is the tableau form of Gessel's expansion of the Schur function in fundamental quasisymmetric functions; see Stanley~\cite[Theorem~7.19.7]{Stanley}.)

\emph{$S(\tau)$ is a subset of $A(s)$.} Equal entries of $\tau$ are standardized from left to right, so if $i$ and $i+1$ receive the same entry, then $i+1$ lies in a column strictly to the right of $i$, that is, $i\in A(s)$.

\emph{The map $\tau\mapsto S(\tau)$ is injective.} The entry of $\tau$ in the box labeled $j$ in $s$ is $1+\#\{i<j:\ i\notin S(\tau)\}$, so $\tau$ is determined by $s$ and $S(\tau)$.

\emph{The map $\tau\mapsto S(\tau)$ is surjective.} For $S\subseteq A(s)$, let $\tau$ be obtained from $s$ by giving the labels $i,i+1$ the same entry for $i\in S$. Since $S\subseteq A(s)$, the labels in a maximal run $i,i+1,\dots,j$ joined by $S$ lie in strictly increasing columns, so they form a horizontal strip. Hence the entries of $\tau$ weakly increase along rows and strictly increase down columns, so $\tau$ is semistandard. Standardizing from left to right returns $s$, so $\std(\tau)=s$ and $S(\tau)=S$.

The labels $i,i+1$ receive the same entry exactly when $i$ is not a partial sum $s_1,\dots,s_{r-1}$ of the content $\alpha$, so by the definition of $w_\alpha$ we have $w_\tau=\prod_{i\in S(\tau)}(q^i-1)$ and
\[
\sum_{\std(\tau)=s}w_\tau=\sum_{S\subseteq A(s)}\ \prod_{i\in S}(q^i-1)=\prod_{i\in A(s)}\bigl(1+(q^i-1)\bigr)=q^{\sum_{i\in A(s)}i}.\qedhere
\]
\end{proof}

\begin{example}
For $\lambda=(2,1)$ the four tableaux of Example~\ref{ex21} have the following standardizations: the first two (contributing $1$ and $q-1$ vectors) both standardize to the tableau $\ytab{1\&2\\3}$, for which $A(s)=\{1\}$ and the total is $q$; the last two (contributing $1$ and $q^2-1$ vectors) both standardize to the tableau $\ytab{1\&3\\2}$, for which $A(s)=\{2\}$ and the total is $q^2$. So the basis splits into sets of sizes $q$ and $q^2$.
\end{example}

\subsection*{The standard basis problem}

Since the sets $\mathcal B_s$ together form a basis of $S^\lambda$, their sizes add up to $\dim S^\lambda$. This sum agrees with the known formulas for the dimension. Transposing the tableaux $s$ gives
\[
\sum_sq^{\maj(s^T)}=\sum_{t\in\mathrm{SYT}(\lambda')}q^{\maj(t)}=q^{n(\lambda')}\,\frac{[n]_q!}{\prod_b[h(b)]_q}
\]
by Macdonald~\cite{Mac}, III.6, Example~2 together with I.5, Example~14, formula~(3). Here $n(\lambda')=\sum_i(i-1)\lambda'_i$, $[m]_q=(q^m-1)/(q-1)$, $[n]_q!=[1]_q\cdots[n]_q$, and $h(b)$ is the hook length of the box $b$ of $\lambda$. The right side is the degree of the unipotent character of $\GL_n(\FF_q)$ attached to $\lambda$~\cite{DJ}, (1.1) (with the partition transposed relative to our convention), which is $\dim S^\lambda$. The middle sum is the decomposition $\dim S=\sum_tq^{f(t)}$ with $f=\maj$~\cite{Thomas} that Dipper and James quote~\cite{DJ}, for the standard tableaux $t=s^T$ of the conjugate shape $\lambda'$ (their Specht module $S^{\lambda'}$ is our $S^\lambda$).

In the formulation of Dipper--Guo~\cite{DG}, the standard basis problem asks for polynomials $g_s$ with $g_s(1)=1$ and, for each standard tableau $s$, a set of $g_s(q)$ vectors, so that the union of these sets is a basis. Read literally, the sets $\mathcal B_s$ solve this, with $g_s(q)=q^{\maj(s^T)}$.

Dipper and James~\cite{DJ} ask for more: the vectors should be naturally attached to the tableau $t=s^T$. They report evidence that there is no natural way of attaching $q^{f(t)}$ vectors to $t$, and they conjecture different weights. For two-part partitions these weights are the rank polynomials of~\cite{BDJL}, whose Theorem~3.15 shows that the leading component of the Specht module attached to a row-standard tableau $t$ has dimension $r(t)$. For example, for $S^{(3,3)}$ (James's labeling) the tableaux $135/246$, $125/346$, $134/256$, $124/356$, $123/456$ have weights $q^3,q^5,q^5,q^7+q^6-q^5,q^9$, instead of the five powers $q^3,q^5,q^6,q^7,q^9$ (Example~2.5 there). These tableaux and polynomials for $S^{(3,3)}$, and for $S^{(2,2)}$ the weights $q^2$ and $q^4$ for $13/24$ and $12/34$, are listed in~\cite[p.~107 and Example~4.10]{DJ} and~\cite[pp.~73--74]{James84}.

What is proved here is therefore only a cardinality statement: a basis partitioned into sets of the sizes $q^{f(t)}$, where the partition is defined combinatorially by standardizing the chain. We do not claim that the sets $\mathcal B_s$ are compatible with the module structure in the sense of Dipper and James, and we would not expect this for these weights. Nor do we claim that the specialization $q\to1$ of the vectors of $\mathcal B_s$ is James's standard polytabloid $e_s$. All that is true is that at $q=1$ every weight $w_\tau$ with a nontrivial strip vanishes, so that only the standard chains $\tau=s$ (all strips single boxes) have nonzero weight.

For two-part partitions the standard basis problem was settled by Brandt, Dipper, James, and Lyle~\cite{BDJL} (see also~\cite{Brandt}). Their basis vectors are Fourier sums over the fillings of a Schubert cell, ordered by the last path in reverse lexicographic order. Lyle~\cite{Lyle} first proved the vanishing of the part attached to a nonstandard tableau, a special case of their dimension theorem. Guo~\cite{Guo,DG} gave a representation-theoretic construction of that standard basis, through the restriction of $S^\lambda$ and $M^\lambda$ to the unitriangular group, together with an explanation of the rank polynomials. James's book~\cite[Ch.~12 (p.~73), Thm.~16.9]{James84} constructs a basis of $S^\lambda$ for arbitrary $\lambda$ by a recursion along the layers of the restriction to the mirabolic subgroup, of the same shape as ours (see the Introduction).

\section{Other fields}\label{secmodular}

In Sections~\ref{secbasis}--\ref{secstd} the coefficient field was $\CC$. In this section we show that the recursive family remains a basis of $S^\lambda$ over every field of characteristic $\ell\neq p$ (including when $\ell$ divides $|G|$), and we relate it to the simple module $D^\lambda$.

\subsection{Independence in every characteristic \texorpdfstring{$\ell\ne p$}{l != p}}\label{subindep}

Let $\KK$ be a field containing a nontrivial $p$th root of unity. Its characteristic $\ell$ is then $0$ or different from $p$ ($\ell$ may divide $|G|$). Let $A=\ZZ[\zeta_p,1/p]$ and choose an $A$-algebra structure on $\KK$. All the objects of the construction ($x_T$, $\Omega_{Y_k}$, $W$, $R_k$, $B(\lambda)$) have coefficients in $A$, and so are defined over $\KK$; the module $M^\lambda_\KK=\KK[X_\lambda]$ carries the symmetric bilinear form $[\,,\,]$ for which the flags are orthonormal. The constants $|N|$, $|Y_k|$, $c$ are powers of $p$, hence units. For a $k$-strip $\lambda/\mu$, define the $\GL_{n-k}$-equivariant map
\[
\Theta_\mu=J^\dagger \Omega_{Y_k}\colon M^\lambda\to M^\mu
\]
where $J^\dagger$ is the adjoint of $J$ from Section~\ref{secdiag}, now taken for the bilinear form $[\,,\,]$ (equivalently, the transpose of the $0/1$ matrix of $J$), and $\Omega_{Y_k}$ acts as an operator on flags. Put $d=\dim_\CC S^\lambda=|B(\lambda)|$, and let $S^\lambda_\KK=\KK\text{-span}\{e_T\}\subseteq M^\lambda_\KK$ be the module of~\cite{AN2} over $\KK$.

\begin{proposition}\label{propids}
For every $k$-strip $\lambda/\mu$ and every $\beta\in A[\GL_{n-k}]$:
\begin{enumerate}
\item[(a)] $\Theta_\mu(\beta \Omega_{Y_k}x_{T'(\lambda,\mu)})=|Y_k|\,c\,\beta x_{T_\mu}$;
\item[(b)] $\Theta_\mu(\beta \Omega_{Y_k}x_{T'(\lambda,\nu)})=0$ for every $k$-strip $\nu$ with $\nu\not\domeq\mu$.
\end{enumerate}
These are identities between vectors with entries in $A$, so they hold over every $\KK$.
\end{proposition}
\begin{proof}
\emph{Reduction to $\CC$.} Both sides of (a) and (b) are vectors in $M^\mu_A$, since $J^\dagger$ (the transpose of a $0/1$ matrix), $\Omega_{Y_k}$, $\beta$, and $x_{T'(\lambda,\nu)}$ all have entries in $A$. Since $A\subset\CC$, an identity between such vectors that holds over $\CC$ holds over $A$, and hence over $\KK$ by the structure map $A\to\KK$. So it suffices to prove (a) and (b) over $\CC$, where the Hermitian form is available and $J^\dagger$ is the adjoint of $J$.

(a) Since $\beta$ commutes with $\Theta_\mu$ and $\Omega_{Y_k}^2=|Y_k|\,\Omega_{Y_k}$ (Section~\ref{subnotation}(i)), we have $\Theta_\mu(\beta \Omega_{Y_k}x_{T'})=|Y_k|\,\beta J^\dagger\Omega_{Y_k}x_{T'}=|Y_k|\,\beta\xi=|Y_k|\,c\,\beta x_{T_\mu}$ by Theorem~\ref{thmD2}.

(b) For $s\in M^\mu$ and $z\in M_k$, using the Hermitian form throughout: $\langle\Theta_\mu z,s\rangle=\langle \Omega_{Y_k}z,Js\rangle=\langle z,\Omega_{Y_k}Js\rangle=\langle z,\Omega_{Y_k}\pi_\lambda Js\rangle$. The second equality holds because $\Omega_{Y_k}$ is self-adjoint over $\CC$ (Section~\ref{subnotation}(v)), and the third because $z\in S^\lambda$. The map $\Omega_{Y_k}\pi_\lambda J\colon M^\mu\to M_k$ is $\GL_{n-k}$-equivariant, so its image lies in the sum of the summands $M_k^\rho$ for which $S^\rho$ is a constituent of $M^\mu$, and these have $\rho\domeq\mu$ (Section~\ref{subflags}(e)). Hence the image of $\Omega_{Y_k}\pi_\lambda J$ lies in $\bigoplus_{\rho\domeq\mu}M_k^\rho$, and $\langle\Theta_\mu z,s\rangle$ depends only on the components $z^\rho$ of $z$ with $\rho\domeq\mu$. Now take $z=\beta \Omega_{Y_k}x_{T'(\lambda,\nu)}$ with $\nu\not\domeq\mu$. Then every component $z^\rho$ with $\rho\domeq\mu$ is zero: a nonzero component has $\rho\trianglerighteq\nu$ (Theorem~\ref{thmT}), and $\rho\domeq\mu$ would then give $\nu\domeq\mu$. Hence $\langle\Theta_\mu z,s\rangle=0$ for all $s$, so $\Theta_\mu z=0$.
\end{proof}

\begin{theorem}[Independence in every characteristic $\ell\ne p$]\label{thmindep}
For $\KK$ as above the vectors $\{\beta x_{T_\lambda}:\beta\in B(\lambda)\}$ are linearly independent in $M^\lambda_\KK$, and form a basis of the submodule $\KK[G]\,x_{T_\lambda}\subseteq M^\lambda_\KK$, which has dimension $d$.
\end{theorem}
\begin{proof}
We use induction on $n$. By Lemma~\ref{lemlayers}, which only uses that $|N_n|$ is invertible, $M^\lambda_\KK=\bigoplus_{k=1}^n\bigoplus_{g\in R_k}g\,\Omega_{Y_k}M^\lambda_\KK$, so it suffices to prove the independence of $\{\beta \Omega_{Y_k}x_{T'(\lambda,\mu)}\}$ ($\mu$ a $k$-strip, $\beta\in B(\mu)$) for each $k$. Take a relation $\sum a_{\nu,\beta}\beta \Omega_{Y_k}x_{T'(\lambda,\nu)}=0$, and let $\mu$ be dominance-minimal among the $\nu$ with some $a_{\nu,\beta}\ne0$. Every other $\nu$ in the support is strictly above $\mu$ or incomparable to $\mu$ (dominance is only a partial order), so $\nu\not\domeq\mu$. Applying $\Theta_\mu$ and using Proposition~\ref{propids} leaves $|Y_k|c\sum_\beta a_{\mu,\beta}\beta x_{T_\mu}=0$. The factor $|Y_k|c$ is a unit, and $\{\beta x_{T_\mu}\}$ is independent by induction, so all $a_{\mu,\beta}=0$, a contradiction.

Finally, consider the matrix whose columns are the flag coordinates of the vectors $g\,x_{T_\lambda}$, $g\in G$; its entries lie in $A$ and, by Theorem~\ref{thmx}, its rank over $\CC$ is $\dim_\CC\CC[G]x_{T_\lambda}=d$. So all its $(d+1)\times(d+1)$ minors vanish in $A$, hence in $\KK$, and $\dim_\KK\KK[G]x_{T_\lambda}\le d=|B(\lambda)|$. The independent family $\{\beta x_{T_\lambda}\}$, which lies in $\KK[G]x_{T_\lambda}$, is therefore a basis of it.
\end{proof}

\begin{remark}[Integral form]\label{remint}
By Theorem~\ref{thmindep}, the family $\{\beta x_{T_\lambda}:\beta\in B(\lambda)\}$ stays linearly independent modulo every maximal ideal $\mathfrak m$ of $A$. Let $L$ be its $A$-span. Independence modulo $\mathfrak m$ means that some $d\times d$ minor of the coordinate matrix of the family is a unit in the localization $A_\mathfrak m$, so $L_\mathfrak m$ is a direct summand of $(M^\lambda_A)_\mathfrak m$. Hence $M^\lambda_A/L$ is torsion free, that is, $L$ is a pure sublattice of $M^\lambda_A$. Since $L\otimes\CC=S^\lambda_\CC$ (Theorem~\ref{thmmain}), we get $L=M^\lambda_A\cap S^\lambda_\CC$. For every $\KK$ as above, the reduction of this lattice in $M^\lambda_\KK$ is then $\KK[G]x_{T_\lambda}$.
\end{remark}

\subsection{The module \texorpdfstring{$S^\lambda_\KK$}{S(lambda,K)} and James's Specht modules}\label{subJames}

Theorem~\ref{thmindep} shows that the family is a basis of $\KK[G]x_T$. We also want it to be a basis of $S^\lambda_\KK=\KK\text{-span}\{e_T\}$, the module of~\cite{AN2}, whose quotient $D^\lambda$ is the simple module. By Corollary~\ref{corkappa}, $k_Tx_T=\kappa e_T$ with $\kappa=|N|$ a unit, so $S^\lambda_\KK\subseteq\KK[G]x_T$. The following are equivalent:
\begin{enumerate}
\item[(i)] $S^\lambda_\KK=\KK[G]x_T$;
\item[(ii)] $x_T\in S^\lambda_\KK$;
\item[(iii)] $\dim_\KK S^\lambda_\KK=d=\dim_\CC S^\lambda$;
\item[(iv)] the $A$-lattice spanned by the $e_T$ is pure at $\ell$.
\end{enumerate}
We refer to these equivalent statements as \textbf{(P)}.

We prove (P) by identifying $S^\lambda_\KK$ with James's Specht module, whose dimension is the same in characteristic $0$ and in characteristic $\ell$. We carry out this identification now.

\medskip
\noindent\textbf{James's Specht modules.} We use the Specht modules of James~\cite{Jam86} for $s=1$ (so his $d=1$ and $k=n$). James works with right modules over a field $F$ in which $q\cdot1_F\ne0$ and which contains a primitive $p$th root of unity; his fields are a number field $K$ and the residue field $\bar K$ of a prime above $\ell$. For a composition $\mu$ of $n$ let $P^+_\mu$ be the block upper triangular matrices with diagonal blocks of sizes $\mu_1,\mu_2,\dots$, let $A(\mu)\subseteq\{1,\dots,n-1\}$ be the set of $i$ such that $i,i+1$ lie in the same row of the row-reading tableau of shape $\mu$, let $U^+$ be the upper unitriangular group, and for $B\subseteq\{1,\dots,n-1\}$ put $\theta_B(u)=\theta\bigl(\sum_{i\in B}u_{i,i+1}\bigr)$, a character of $U^+$. Write $\theta_\mu=\theta_{A(\mu)}$ and $E(\mu)=|U^+|^{-1}\sum_{u\in U^+}\theta_\mu(u^{-1})u$~\cite[Def.~1.3, (1.6)]{Jam86}. James's module $M(1,\mu)$ is the permutation module on the right cosets $P^+_\mu\backslash G$~\cite[\S2]{Jam86}, and his Specht module is
\[
S(1,\mu)=M(1,\mu)\,E(\mu')\,\KK G,\qquad \mu'=\text{the conjugate of }\mu,
\]
the right submodule generated by $M(1,\mu)E(\mu')$~\cite[Def.~7.7]{Jam86}. For a composition $\mu$, the conjugate $\mu'$ is defined by the same formula as for partitions, so it is the conjugate of the partition obtained by sorting the parts of $\mu$. We use the following results of~\cite{Jam86} for $s=d=1$.
\begin{itemize}
\item[(J1)] [Thm.~7.1(ii)] $M(1,\mu)E(\mu')$ is one dimensional, for every composition $\mu$ of $n$.
\item[(J2)] [(1.4)] If $\alpha'=\mu'$ (that is, $\alpha$ is obtained from $\mu$ by rearranging the parts) then $\dim ME(\alpha)=\dim ME(\mu)$ for every $\KK G$-module $M$.
\item[(J3)] [Thm.~7.19(i)] For a partition $\mu$ of $n$, $\dim S_F(1,\mu)$ is the same for $F=K$ and $F=\bar K$.
\end{itemize}
Of these, only (J3) is used without proof: (J1) is reproved in Remark~\ref{remJ1}, and (J2) in Lemma~\ref{lemJ2}. The choice of the nontrivial character $\theta$ plays no role~\cite[pp.~239--240, before (1.9)]{Jam86}.

\begin{remark}[On the sentence before (1.4) in~\cite{Jam86}]\label{remJ2}
The sentence before (1.4) in~\cite{Jam86} asserts that, for $\alpha'=\mu'$, there is $g\in G$ with $g^{-1}E(\mu)g=E(\alpha)$. Taken literally this cannot hold. Take $q=2$, $n=3$, $\mu=(2,1)$, $\alpha=(1,2)$, over $\CC$: then $E(\mu)$ and $E(\alpha)$ have every coefficient on $U^+$ nonzero, so a conjugating element $g$ would have to normalize $U^+$, which for $q=2$ forces $g\in U^+$; but conjugation by $U^+$ fixes every linear character of $U^+$, while $\theta_\mu(u)=\theta(u_{12})\neq\theta(u_{23})=\theta_\alpha(u)$ for the element $u$ with $u_{12}=1$ and $u_{23}=0$. We therefore use only the numerical statement (1.4), i.e., (J2), and prove (J2) and the consequence for generated submodules that we need in Lemma~\ref{lemJ2}.
\end{remark}

\begin{remark}[Related results in the literature]
The facts (J1) and (J3), the unique maximal submodule of $S(1,\mu)$, and the Kostka-number filtration of $M(1,\mu)$ are restated in~\cite{DJid}, (4.1)--(4.3).

James's earlier construction of unipotent modules~\cite{James82}, \S3, realizes the permutation module on the cosets of the block upper triangular parabolic subgroup as a left ideal $(\FF G)\bar P_J$ of the group algebra. There the analogue of (J3), that the module defined over a prime field of characteristic $\ell\ne p$ is the reduction of the one defined over $\mathbb Q$, is only conjectured (end of \S4). It is proved for the modules $S(1,\mu)$ of~\cite{Jam86}.

For general $\ell\nmid q$, Theorem~3.5e and Remark~3.5f of~\cite{BDK} identify James's module $S(\sigma,\lambda)$ with a standard module $\Delta(\sigma,\lambda)$ of a quasi-hereditary quotient of $FG_n$. Its composition multiplicities are those of Weyl modules of the $q$-Schur algebra, $[\Delta(\sigma,\lambda):L(\sigma,\mu)]=[\Delta_h(\lambda'):L_h(\mu')]$ (\cite{BDK}, (3.5b)(ii)). It has simple head $L(\sigma,\lambda)$, and all its other composition factors are $L(\sigma,\mu)$ with $\mu>\lambda$.
\end{remark}

\begin{lemma}\label{lemJ2}
Let $\alpha,\mu$ be compositions of $n$ with $\alpha'=\mu'$. Then for every finite-dimensional right $\KK G$-module $M$ (with $\KK$ as above):
\begin{enumerate}
\item[(a)] $\dim ME(\alpha)=\dim ME(\mu)$ (this is (J2));
\item[(b)] $ME(\alpha)\,\KK G=ME(\mu)\,\KK G$.
\end{enumerate}
\end{lemma}
\begin{proof}
(a) The dimensions of the images of the idempotents are unchanged by extension of scalars. So we may extend $\KK$ to a splitting field for $G$, and choose characteristic-$0$ lifts of the roots of unity that are compatible with the chosen lift of $\theta$. Let $\tilde\theta_\alpha$ be the complex lift of $\theta_\alpha$ and $\Psi_\alpha=\Ind_{U^+}^{G}\tilde\theta_\alpha$. Write $U^+=N_\alpha U_{L_\alpha}$ (the unipotent radical of the parabolic subgroup of type $\alpha$ and the unitriangular group of its Levi factor $L_\alpha$). Since $\theta_\alpha$ is trivial on $N_\alpha$ and generic on each unitriangular diagonal block, transitivity of induction gives $\Psi_\alpha=R_{L_\alpha}^{G}(\Gamma_{\alpha_1}\otimes\cdots\otimes\Gamma_{\alpha_r})$, where $\Gamma_a$ is the Gelfand--Graev character of $\GL_a$ and $R$ is Harish-Chandra induction. Harish-Chandra multiplication on the full complex representation ring is commutative~\cite{Z}, 9.1, so $\Psi_\alpha=\Psi_\mu$. Now let $\varphi_M$ be the character of $M$ if $\mathrm{char}\,\KK=0$ (the ordinary character, with the chosen compatible identification of roots of unity), and the Brauer character of $M$ if $\mathrm{char}\,\KK=\ell>0$, lifted to characteristic $0$ with the same roots of unity and extended by $0$ to the $\ell$-singular elements (nothing needs to be extended in characteristic $0$). If $\ell>0$, every element of $U^+$ has order prime to $\ell$, so $ME(\alpha)$ is the $\theta_\alpha$-eigenspace of $U^+$ in $M$ and
\[
\dim ME(\alpha)=\frac1{|U^+|}\sum_{u\in U^+}\tilde\theta_\alpha(u^{-1})\varphi_M(u)=\langle\varphi_M,\Psi_\alpha\rangle_G
\]
(Frobenius reciprocity; $\Psi_\alpha$ is supported on the conjugates of $U^+$, which are $\ell$-regular, so extending $\varphi_M$ by $0$ is harmless). The right side depends only on $\Psi_\alpha=\Psi_\mu$.

(b) Put $R=\KK G$, $e=E(\alpha)$, $f=E(\mu)$, idempotents of $R$. Apply (a) to the module $M=R/ReR$: $Me=0$, so $Mf=0$, i.e., the class of $f$ in $R/ReR$ vanishes ($Mf$ contains the class of $1\cdot f$). So $f\in ReR$, and symmetrically $e\in RfR$. Hence $ReR=RfR$, and for every module $M$, $MeR=M(ReR)=M(RfR)=MfR$.
\end{proof}

\begin{theorem}[(P); identification with James's module]\label{thmP}
Let $\KK$ be as above, $T=T_\lambda$ the row-reading tableau, $y=\sum_{n\in N(T)}n\,e_T$, and let $\mathcal X=X_\lambda$ be the set of flags in $\FF_q^n$ of type $\lambda'$ (Section~\ref{subflags}), so that $M^\lambda_\KK=\KK[\mathcal X]$. Then:
\begin{enumerate}
\item[(a)] $y\neq0$, and $\{z\in M^\lambda_\KK:\ vz=\theta_\lambda(v)z\ \ \forall v\in U^+\}=\KK y$ is one dimensional.
\item[(b)] $\KK[G]\,y$ is James's Specht module $S(1,\lambda')$, where the right module $S(1,\lambda')$ is regarded as a left module through the antiautomorphism $\sigma(g)=w_0g^{T}w_0$ of $G$ ($w_0$ the permutation matrix reversing $1,\dots,n$), which reflects a matrix across its antidiagonal.
\item[(c)] $S^\lambda_\KK=\KK[G]e_T=\KK[G]\,y=\KK[G]x_T$ and $\dim S^\lambda_\KK=d$; that is, (P) holds.
\item[(d)] $\mathrm{Rad}=S^\lambda_\KK\cap(S^\lambda_\KK)^\perp$ is the unique maximal submodule of $S^\lambda_\KK$ (\cite{AN2}, Cor.~2); by (b) it is James's $S(1,\lambda')^{\max}$, so $D^\lambda=S^\lambda_\KK/\mathrm{Rad}$ is simple and is James's $D(1,\lambda')$.
\end{enumerate}
\end{theorem}

\begin{proof}
\emph{Step 1: translation to left modules.} For a composition $\mu=(\mu_1,\dots,\mu_r)$ write $\rev\mu=(\mu_r,\dots,\mu_1)$, and put $\nu=(\lambda'_{\lambda_1},\dots,\lambda'_1)=\rev\lambda'$. As an increasing flag, $F_T$ has type $\nu$, so $\mathcal X\cong G/P^+_\nu$ as $G$-sets (the standard flag $\mathrm{span}\{e_1,\dots,e_{\nu_1}\}\subset\mathrm{span}\{e_1,\dots,e_{\nu_1+\nu_2}\}\subset\cdots$ has stabilizer $P^+_\nu$). Write $W=\KK[G/P^+_\nu]\cong\KK[\mathcal X]$, with $G$ acting by left multiplication. The antiautomorphism $\sigma$ of~(b) satisfies $\sigma(U^+)=U^+$ and $\sigma(P^+_{\lambda'})=P^+_\nu$, and $\sigma(u)_{k,k+1}=u_{n-k,\,n-k+1}$. The map $P^+_{\lambda'}g\mapsto\sigma(g)P^+_\nu$ is a bijection $P^+_{\lambda'}\backslash G\to G/P^+_\nu$ with $\sigma(h)\cdot(\text{image of }P^+_{\lambda'}g)=\text{image of }P^+_{\lambda'}gh$. So $M(1,\lambda')\cong W$ with the right action of $h$ becoming the left action of $\sigma(h)$, and a right submodule generated by a vector $z$ becomes $\KK[G]z$. For a composition $\mu$, the elements $m$ of $M(1,\lambda')E(\mu)$ are those with $m\cdot u=\theta_\mu(u)m$ for all $u\in U^+$ (a direct check from the definition of $E(\mu)$), and the corresponding vectors $z$ are those with $\sigma(u)z=\theta_\mu(u)z$ for all $u\in U^+$; since $\theta_\mu(u)=\theta_{\,n-A(\mu)}(\sigma(u))$ and $n-A(\mu)=A(\rev\mu)$, we get
\[
M(1,\lambda')E(\mu)\ \cong\ W_{\rev\mu}:=\{z\in W:\ vz=\theta_{\rev\mu}(v)z\ \ \forall v\in U^+\}.
\]
Take $\mu=\rev\lambda$. Then $\mu'=\lambda'$, so by (J1) and Lemma~\ref{lemJ2}(a) (applied to $\alpha=\lambda$ and $\mu=\rev\lambda$, which have the same conjugate $\lambda'$) we have $\dim W_\lambda=\dim M(1,\lambda')E(\rev\lambda)=\dim M(1,\lambda')E(\lambda)=1$, and the submodule of $W$ generated by $W_\lambda$ is James's $S(1,\lambda')=M(1,\lambda')E(\lambda)\KK G$ (equal to $M(1,\lambda')E(\rev\lambda)\KK G$ by Lemma~\ref{lemJ2}(b)).

\emph{Step 2: $y\in W_\lambda$.} We have $U^+=N\,U_R$, where $N=N(T)$ is the unipotent radical of the block upper triangular parabolic subgroup with the rows as blocks and $U_R$ is the unitriangular group of its Levi factor; $N$ is normal in $U^+$ and $N\cap U_R=1$. Moreover $U_R\subseteq U(T)$ and the pairs in $X(T)$ are, in terms of labels, the pairs $(i,i+1)$ with $i\in A(\lambda)$, so $\psi_T(u)=\theta_\lambda(u)$ for $u\in U(T)\cap U^+$, and $ue_T=\psi_T(u)e_T$ for $u\in U(T)$ by (A1). For $r\in U_R$: $r\,y=\sum_{n\in N}(rnr^{-1})\,r\,e_T=\psi_T(r)\,y$. For $n'\in N$: $n'y=y$, and $\theta_\lambda(n')=1$ since $n'_{i,i+1}=0$ for $i\in A(\lambda)$. Hence for $v=n'r$: $vy=\psi_T(r)y=\theta_\lambda(v)y$.

\emph{Step 3: $y\ne0$.} The coefficient of $[F_T]$ in $y=\sum_{n\in N}\sum_{u\in U(T)}\psi_T(u)^{-1}[nuF_T]$ is $\sum\psi_T(u)^{-1}$ over the pairs $(n,u)\in N\times U(T)$ with $nuF_T=F_T$. For such a pair, $uF_T=n^{-1}F_T\in NF_T$, so $u\in A=N\cap U(T)$ by Lemma~\ref{lemorbits}(b). Then $nu\in N$ fixes $F_T$, so $nu\in N\cap P(T)$. Conversely, each $(a,p)\in A\times(N\cap P)$ gives such a pair, with $u=a$ and $n=pa^{-1}$. So these pairs correspond bijectively to $A\times(N\cap P)$. As $\psi_T=1$ on $A$ (no entries inside a row), the coefficient is $|A|\,|N\cap P|=|N|=q^{\sum_{r<s}\lambda_r\lambda_s}$ (Lemma~\ref{lemorbits}(a)), a power of $p$, which is nonzero in $\KK$. By Step~1, $W_\lambda$ is one dimensional, and it contains $y$ by Step~2. Since $y\ne0$, we get $W_\lambda=\KK y$, which is (a), and the submodule generated by $W_\lambda$ is $\KK[G]y=S(1,\lambda')$, which is (b).

\emph{Step 4: dimensions; proof of (c).} First, the dimensions are unchanged under extension of the coefficient field within a fixed characteristic. The module $W$, the idempotents, and $y$ are all defined over the subfield $k_0\subseteq\KK$ generated by the values of $\theta$, and the dimensions are ranks of matrices over $k_0$, so they do not change under extension of scalars. Identify $\theta$ with James's character through an embedding of $k_0$ (the choice does not matter, as noted after (J3)). Then these dimensions agree with those of James's modules over his fields $K$ (of characteristic $0$) and $\bar K$. The comparison between the two characteristics is (J3).

Next we compute $d$ in characteristic $0$. The module $\CC[G]y$ is a nonzero submodule of $S^\lambda_\CC=\CC[G]e_T$, which is irreducible by (A3). So $S^\lambda_\CC=\CC[G]y$, which is $S(1,\lambda')_K$ after extension from $K$ to $\CC$, and $d=\dim S(1,\lambda')_K$. By (J3), $\dim S(1,\lambda')_{\bar K}=d$ as well. Hence $\dim_\KK\KK[G]y=d$ for every $\KK$ as above.

Now $\KK[G]y\subseteq\KK[G]e_T=S^\lambda_\KK$, because $y\in\KK[G]e_T$. Also $\dim S^\lambda_\KK\le d$: the matrix whose columns are the flag coordinates of the vectors $g\,e_T$ has entries in $A$ and rank $d$ over $\CC$, so its $(d+1)\times(d+1)$ minors vanish in $A$, hence in $\KK$. Therefore $S^\lambda_\KK=\KK[G]y$, and this has dimension $d$.

Finally, $e_T=\kappa^{-1}k_Tx_T$ gives $S^\lambda_\KK\subseteq\KK[G]x_T$, and $\KK[G]x_T$ has dimension $d$ by Theorem~\ref{thmindep}. So $S^\lambda_\KK=\KK[G]x_T$, which completes (c).

\emph{Step 5: proof of (d).} Let $\bar T$ be the tableau obtained from $T$ by replacing each entry in an odd column by its negative. By~\cite{AN2}, Lemma~6, $\bar T$ has the same flag as $T$ (so $m_{\bar T}=m_T$), $U(\bar T)=U(T)$ and $\psi_{\bar T}(u)=\psi_T(u^{-1})$; hence $e_{\bar T}=\sum_{u\in U(T)}\psi_T(u)\,[uF_T]$. For every $x\in M^\lambda_\KK$ we claim
\[
k_Tx=[x,e_{\bar T}]\;e_T .
\]
Indeed $k_Tx\in\KK e_T$ by (A2), the coefficient of $[F_T]$ in $e_T$ is $1$, and the coefficient of $[F_T]$ in $k_Tx=\sum_u\psi_T(u^{-1})\,ux$ is $\sum_u\psi_T(u^{-1})x(u^{-1}F_T)=\sum_{v\in U(T)}\psi_T(v)\,x(vF_T)=[x,e_{\bar T}]$. Now let $V$ be a submodule of $M^\lambda_\KK$. If $k_Tx\ne0$ for some $x\in V$ and some $T$, then $e_T\in V$, so $S^\lambda_\KK\subseteq V$. Otherwise $[x,e_{\bar T}]=0$ for all $x\in V$ and all $T$; since $T\mapsto\bar T$ is a bijection on the set of tableaux of shape $\lambda$, this says $V\subseteq(S^\lambda_\KK)^\perp$. Moreover $S^\lambda_\KK\not\subseteq(S^\lambda_\KK)^\perp$: since the flags $uF_T$ are distinct, $[e_T,e_{\bar T}]=\sum_{u\in U(T)}\psi_T(u^{-1})\psi_T(u)=|U(T)|\ne0$, and $e_T,e_{\bar T}\in S^\lambda_\KK$. So every proper submodule of $S^\lambda_\KK$ lies in $\mathrm{Rad}$, and $\mathrm{Rad}\ne S^\lambda_\KK$. This is the unique maximal submodule; by (b) it coincides with James's (Thm.~7.8 of~\cite{Jam86}).
\end{proof}

Thus $S^\lambda_\KK=\KK[G]x_T$ for every $\KK$ as above, and $D^\lambda$ (the module of~\cite{AN2}) is James's simple module $D(1,\lambda')$. The irreducibility of $D^\lambda$ and the fact that $\mathrm{Rad}$ is the unique maximal submodule are Corollary~2 of~\cite{AN2}; what Theorem~\ref{thmP} adds is that $S^\lambda_\KK$ is James's Specht module and has the characteristic-$0$ dimension, which~\cite{AN2} does not prove (it identifies the simple modules $D^\lambda\cong D(1,\lambda')$ only, through James's characterization of them).

\begin{remark}[On (J1)]\label{remJ1}
The one-dimensionality of $W_\lambda$ (James's (J1) for $\mu=\rev\lambda$) can be checked without~\cite{Jam86}. Since $|U^+|$ is invertible in $\KK$ and $\KK$ contains a primitive $p$th root of unity, the dimension of the $\theta_\lambda$-eigenspace in the permutation module $\KK[\mathcal X]$ is the number of $U^+$-orbits on $\mathcal X$ whose stabilizers lie in $\ker\theta_\lambda$, which does not depend on $\KK$. Over $\CC$ it is $\langle\Ind_{P_{\lambda'}}^G1,\Psi^\lambda\rangle$, where $\Psi^\lambda=\Ind_{U^+}^G\theta_\lambda$ is the degenerate Gelfand--Graev character of Zelevinsky, and by (H1) and $\langle\Psi^\lambda,\chi^\kappa\rangle=K_{\kappa\lambda}$~(\cite{Z}, \S12.1; \cite{AN2}, Prop.~3) this is $\sum_\kappa K_{\kappa'\lambda'}K_{\kappa\lambda}=1$, because $K_{\kappa\lambda}\ne0$ forces $\kappa\trianglerighteq\lambda$ and $K_{\kappa'\lambda'}\ne0$ forces $\kappa\domeq\lambda$.
\end{remark}

\subsection{Bases of the simple modules \texorpdfstring{$D^\lambda$}{D(lambda)}}\label{subD}

Let $S=S^\lambda_\KK$, $\mathrm{Rad}=S\cap S^\perp$ (a $G$-submodule since the form is $G$-invariant) and $D^\lambda=S/\mathrm{Rad}$. If $S$ is semisimple (for instance if $\ell\nmid|G|$) then $\mathrm{Rad}=0$ and $D^\lambda=S^\lambda$.

\begin{lemma}\label{lemrad}
Let $M_k=\Omega_{Y_k}S$ and let $V\subseteq S$ be any $\KK[P_n]$-submodule. Then $\Omega_{Y_k}V=V\cap M_k$ and
\[
V=\bigoplus_{k=1}^n\ \bigoplus_{g\in R_k}g\,(V\cap M_k).
\]
In particular, for $V=\mathrm{Rad}$ the quotient $D^\lambda=S/\mathrm{Rad}$ satisfies
\[
D^\lambda=\bigoplus_{k=1}^n\ \bigoplus_{g\in R_k}g\,D_k,\qquad D_k:=\Omega_{Y_k}D^\lambda=M_k/(M_k\cap\mathrm{Rad}),
\]
where $D_k$ is a $\GL_{n-k}$-module, and $\dim D^\lambda=\sum_{k=1}^n|R_k|\dim D_k$.
\end{lemma}
\begin{proof}
Lemma~\ref{lemlayers} holds for every $\KK[P_n]$-module, in particular for $V$ and for $S/V$. If $m\in V\cap M_k$ then $m=|Y_k|^{-1}\Omega_{Y_k}m\in\Omega_{Y_k}V$, and the other inclusion is clear. For the last assertion take $V=\mathrm{Rad}$, which is a $G$-submodule. Since $\Omega_{Y_k}$ sends the class of $s\in S$ to the class of $\Omega_{Y_k}s$, the module $\Omega_{Y_k}D^\lambda$ consists of the classes $m+\mathrm{Rad}$ with $m\in M_k$. Two such classes $m+\mathrm{Rad}$ and $m'+\mathrm{Rad}$ agree exactly when $m-m'\in M_k\cap\mathrm{Rad}$, so $\Omega_{Y_k}D^\lambda=M_k/(M_k\cap\mathrm{Rad})$.
\end{proof}

\begin{proposition}[Bases of $D^\lambda$ inside the recursive family]\label{propD}
For $v\in S$ write $\bar v=v+\mathrm{Rad}$ for its class in $D^\lambda$. For each $1\le k\le n$ choose a set $I_k$ of pairs $(\mu,\beta)$ ($\lambda/\mu$ a $k$-strip, $\beta\in B(\mu)$) such that $\{\overline{\beta \Omega_{Y_k}x_{T'(\lambda,\mu)}}\}_{(\mu,\beta)\in I_k}$ is a basis of $D_k$ (such sets exist because, by Theorem~\ref{thmP}, the family spans $M_k$). Then $\{\overline{g\,\beta \Omega_{Y_k}x_{T'(\lambda,\mu)}}:\ 1\le k\le n,\ (\mu,\beta)\in I_k,\ g\in R_k\}$ is a basis of $D^\lambda$.
\end{proposition}
\begin{proof}
By Theorems~\ref{thmindep} and~\ref{thmP} the vectors $g\beta \Omega_{Y_k}x_{T'(\lambda,\mu)}$ ($1\le k\le n$, $\mu$ a $k$-strip, $\beta\in B(\mu)$, $g\in R_k$) form a basis of $S=\KK[G]x_T$. Each of them lies in one summand $g\,M_k$ of the decomposition $S=\bigoplus_{k=1}^n\bigoplus_{g\in R_k}g\,M_k$ of Lemma~\ref{lemlayers}. A basis of a direct sum whose elements each lie in a single summand restricts to a basis of each summand, so for each $k$ the vectors $\beta \Omega_{Y_k}x_{T'(\lambda,\mu)}$ form a basis of $M_k$. Since $D_k$ is a quotient of $M_k$ (Lemma~\ref{lemrad}), the sets $I_k$ exist. Because $\mathrm{Rad}$ is a $G$-submodule, $\overline{gv}=g\bar v$, so translating the basis of $D_k$ by $R_k$ gives a basis of $g\,D_k$. By Lemma~\ref{lemrad}, $D^\lambda=\bigoplus_{k=1}^n\bigoplus_{g\in R_k}g\,D_k$, and the proposition follows.
\end{proof}

So the images in $D^\lambda$ of a subfamily of the recursive basis of $S^\lambda$, chosen layer by layer, form a basis of $D^\lambda$. Which subfamily is determined by the radical $\mathrm{Rad}=S\cap S^\perp$, that is, by the bilinear form $[\,,\,]$ on $S$; nothing in the construction determines it. Write $b_1,\dots,b_d$ for the recursive basis of $S$, and call the matrix $([b_i,b_j])$ its \emph{Gram matrix}. The radical consists of the vectors $\sum_ic_ib_i$ whose coefficient vector lies in the kernel of the Gram matrix, so $\dim D^\lambda=d-\dim\mathrm{Rad}$ is the rank of the Gram matrix over $\KK$.

\begin{proposition}[Specht filtration of the layers]\label{propfilt}
Fix $k$, and let $\mathcal S$ be the set of partitions $\mu$ such that $\lambda/\mu$ is a $k$-strip, partially ordered by dominance. For $\mathcal U\subseteq\mathcal S$ let
\[
F_{\mathcal U}=\operatorname{span}\{\beta\,\Omega_{Y_k}x_{T'(\lambda,\nu)}:\ \nu\in\mathcal U,\ \beta\in B(\nu)\}\subseteq M_k,
\]
and for $\mu\in\mathcal S$ put $F_{\trianglerighteq\mu}=F_{\mathcal U}$ with $\mathcal U=\{\nu\in\mathcal S:\nu\trianglerighteq\mu\}$, and $F_{\triangleright\mu}=F_{\mathcal U}$ with $\mathcal U=\{\nu\in\mathcal S:\nu\triangleright\mu\}$. Then:
\begin{enumerate}
\item[(a)] If $\mathcal U$ is an upper set of $\mathcal S$ (that is, $\mu\in\mathcal U$ and $\nu\in\mathcal S$, $\nu\trianglerighteq\mu$ imply $\nu\in\mathcal U$), then $F_{\mathcal U}$ is a $\KK[\GL_{n-k}]$-submodule of $M_k$.
\item[(b)] For every $\mu\in\mathcal S$, $\Theta_\mu$ induces an isomorphism $F_{\trianglerighteq\mu}/F_{\triangleright\mu}\cong\KK[\GL_{n-k}]x_{T_\mu}$ (which is $S^\mu_\KK$ by Theorem~\ref{thmP}).
\end{enumerate}
\end{proposition}
\begin{proof}
(a) Over $\CC$, Theorems~\ref{thmT} and~\ref{thmD2} show that the span of these vectors is $\bigoplus_{\nu\in\mathcal U}M_k^\nu$, which is $\GL_{n-k}$-stable. The vectors are independent over every residue field of $A$ (the proof of Theorem~\ref{thmindep} establishes the independence of the layer family itself, before translating by $R_k$), so as in Remark~\ref{remint} their $A$-span is pure. A pure lattice whose $\CC$-span is $\GL_{n-k}$-stable is itself $A[\GL_{n-k}]$-stable, hence so is its reduction $F_{\mathcal U}$.

(b) By Proposition~\ref{propids}(b), $\Theta_\mu$ kills $F_{\triangleright\mu}$, and by Proposition~\ref{propids}(a) it sends $\beta\,\Omega_{Y_k}x_{T'(\lambda,\mu)}$ to a unit multiple of $\beta x_{T_\mu}$. The generating vectors are independent (Theorem~\ref{thmindep}), so $F_{\trianglerighteq\mu}$ has a basis consisting of a basis of $F_{\triangleright\mu}$ together with the vectors $\beta\,\Omega_{Y_k}x_{T'(\lambda,\mu)}$, $\beta\in B(\mu)$. Their images form a basis of $\KK[\GL_{n-k}]x_{T_\mu}$ (Theorem~\ref{thmindep} for $\mu$), so the induced map takes a basis to a basis. It is $\GL_{n-k}$-equivariant because $\Theta_\mu$ is.
\end{proof}

So $\Omega_{Y_k}S^\lambda_\KK$ has a filtration with sections the modules $\KK[\GL_{n-k}]x_{T_\mu}$ (these are the Specht modules $S^\mu_\KK$, by Theorem~\ref{thmP}), $\mu$ running over the $k$-strips; this is the modular form of Proposition~\ref{proplayer}. Compare \cite[Thm.~16.9]{James84}, which gives a series with sections $S_{\mu_R}$ for James's Specht modules over every field of characteristic not dividing $q$. The splitting into isotypic summands, the isomorphism $\Phi_{\rm can}$, and the Hermitian-form arguments of the earlier sections are not available when $\ell\mid|G|$, which is why Theorem~\ref{thmindep} is proved with $\Theta_\mu$ instead.

\subsection{The symmetric group}\label{subsym}

For the symmetric group the situation is parallel. James's standard polytabloids form a basis of the Specht module $S^\lambda$ in every characteristic~\cite{James78}, and $D^\lambda\ne0$ exactly for $\ell$-regular $\lambda$, but the simple modules $D^\lambda=S^\lambda/(S^\lambda\cap(S^\lambda)^\perp)$ ($\lambda$ $\ell$-regular, $\ell$ the characteristic) have no known uniform basis: $\dim D^\lambda$ is the rank of the Gram matrix mod $\ell$ (the number of elementary divisors of the Gram matrix not divisible by $\ell$), which has been computed in tables (L\"ubeck, $n<13$ and further) and determined in closed form only for special families~\cite{Lubeck,KM,KN}. For the Hecke algebra analogue, Gram matrices over $\ZZ[v,v^{-1}]$ (which is not a principal ideal domain) do not in general have elementary divisors~\cite{KM}. Adjustment matrices compare the decomposition numbers in characteristic $0$ and in characteristic $\ell$. James conjectured that they are trivial for blocks of weight less than $\ell$. The conjecture is known for weights $3$ and $4$ when $\ell\ge5$~\cite{Fayers,FayersFour,Low}, and it is false in general~\cite{Williamson}.

Dipper and James~\cite{DJ} observe that the unipotent Specht modules of $\GL_n$ carry information about the irreducible representations of symmetric groups in arbitrary characteristic. The unipotent decomposition matrix of $\GL_n$ contains that of the $q$-Schur algebra (Dipper--James~\cite{DJ89}; see also~\cite{BDK} for the relation with quantum $\GL_n$), and hence that of the Hecke algebra $H_q(S_n)$. For $q\equiv1\pmod\ell$, the Hecke algebra $H_q(S_n)$ is the group algebra $\KK S_n$. Here $e$ denotes the quantum characteristic, the least $r\ge2$ with $1+q+\cdots+q^{r-1}=0$ in $\KK$ ($e=\infty$ if there is none). If $q\not\equiv1\pmod\ell$, then $e$ is the multiplicative order of $q$ modulo $\ell$. If $q\equiv1\pmod\ell$, which is the symmetric group case, then $e=\ell$.

Now delete from the unipotent block of the decomposition matrix of $\GL_n$ the columns indexed by the $D^\mu$ with $\mu'$ $e$-singular, and relabel the rows by $\lambda\mapsto\lambda'$ and the columns by $\mu\mapsto\mu'$. The result is the decomposition matrix of the Hecke algebra, with entries $d^H_{\alpha\beta}=[S^{\alpha'}:D^{\beta'}]$ for $\beta$ $e$-regular. See~\cite{Jam86} (Thm.~8.2(iii)) and~\cite{DJid} (\S6); for the Hecke algebra itself see~\cite[Thms.~7.6, 7.7]{DJH} (James's labels are the transposes of ours).

The recursive basis has exactly this structure. Proposition~\ref{propD} says that, in nondefining characteristic, $D^\lambda$ of $\GL_n$ has a basis consisting of the images of a layerwise subfamily of the basis of $S^\lambda$, with the subfamily dictated by the Gram radical, just as for $S_n$. The bridge to $S_n$ is the following.

\begin{proposition}[The $U$-fixed vectors]\label{propU}
For a subgroup $H\le G$ and a $\KK[G]$-module $X$, write $X^H$ for the space of $H$-fixed vectors. Let $U=U_n$ be the upper unitriangular group and $B$ the Borel subgroup of upper triangular matrices. Let $\tau$ range over the standard chains of $\lambda$ (all strips single boxes, so $R_1=\{1\}$ and $g=1$); there are $f^\lambda=\#\mathrm{SYT}(\lambda)$ of them. Then:
\begin{enumerate}
\item[(a)] Each $x_\tau$ is fixed by $B$, in particular by $U$.
\item[(b)] The vectors $x_\tau$ form a basis of $(\KK[G]x_T)^U$, which has dimension $f^\lambda$ for every $\KK$ as above. By Theorem~\ref{thmP}, this is $(S^\lambda_\KK)^U$, and $(S^\lambda_\KK)^U=(S^\lambda_\KK)^B$.
\end{enumerate}
\end{proposition}
\begin{proof}
(a) Let $\tau$ be a standard chain. Then
\[
x_\tau=E_{n-1}W_{n-1}\cdots E_1W_1\,E_0W_0\,x_{T_\lambda},
\]
where $W_j\in\GL_{n-j}$ and $E_j=\sum_{u\in N_{n-j}}u$ ($N_m$ as in Section~\ref{seclayers}), which is $|N_{n-j}|$ times the projection onto $N_{n-j}$-fixed vectors. The factors act from right to left. Fix $j$. The vector $E_j(W_j\cdots E_0W_0\,x_{T_\lambda})$ is $N_{n-j}$-fixed, and the factors $E_{j+1}W_{j+1},\dots,E_{n-1}W_{n-1}$ that act afterwards lie in $\GL_{n-j-1}$, which normalizes $N_{n-j}$. So $x_\tau$ is $N_{n-j}$-fixed for every $j$. Since $U=N_nN_{n-1}\cdots N_2$, the vector $x_\tau$ is $U$-fixed.

It remains to see that $x_\tau$ is fixed by the diagonal torus $T_n$. Let $z=\sum_{n\in N(T_\lambda)}n\,m_{T_\lambda}$. Every diagonal matrix $t$ fixes $z$ (it fixes the coordinate flag $F_{T_\lambda}$ and normalizes $N(T_\lambda)$). For a permutation matrix $r\in R(T_\lambda)$ the matrix $r^{-1}tr$ is again diagonal, so $t(rz)=r(r^{-1}tr)z=rz$. Hence $x_{T_\lambda}=\sum_r\sgn(r)\,rz$ is fixed by $T_n$. Each $E_j$ is invariant under conjugation by $T_n$, and each $W_j$ normalizes $T_n$, so $x_\tau$ is $T_n$-fixed. (This concerns the standard chains only; for general chains the Fourier sums involve nontrivial characters.)

(b) The $f^\lambda$ vectors $x_\tau$ are independent (Theorem~\ref{thmindep}), so it suffices to show $\dim(\KK[G]x_T)^U=f^\lambda$. We use induction on $n$; the case $n\le1$ is trivial. Since $|N_n|$ is invertible, $(\KK[G]x_T)^{N_n}=\Omega_{Y_1}\KK[G]x_T=M_1$, and since $U=N_nU_{n-1}$, the space $(\KK[G]x_T)^U$ is the space of $U_{n-1}$-fixed vectors of $M_1$. By Proposition~\ref{propfilt}, $M_1$ has a filtration by $\KK[\GL_{n-1}]$-modules whose sections are the modules $\KK[\GL_{n-1}]x_{T_\mu}$, $\mu=\lambda-\square$. Taking $U_{n-1}$-fixed vectors preserves exact sequences, because $|U_{n-1}|$ is invertible in $\KK$. Therefore $\dim(\KK[G]x_T)^U=\sum_\mu\dim(\KK[\GL_{n-1}]x_{T_\mu})^{U_{n-1}}=\sum_\mu f^\mu=f^\lambda$ by the inductive hypothesis. Hence the $x_\tau$ form a basis of $(\KK[G]x_T)^U$. Finally $\KK[G]x_T=S^\lambda_\KK$ by Theorem~\ref{thmP}, and since the $x_\tau$ are $B$-fixed by (a), $(S^\lambda_\KK)^U=(S^\lambda_\KK)^B$.
\end{proof}

\begin{proposition}[Gram rank of the $U$-fixed vectors]\label{propGramU}
Let $S=S^\lambda_\KK$, $\mathrm{Rad}=S\cap S^\perp$, and $D^\lambda=S/\mathrm{Rad}$. The Gram matrix $([x_\tau,x_{\tau'}])$ of the vectors $x_\tau$, $\tau$ a standard chain, with respect to the symmetric bilinear form has rank $\dim(D^\lambda)^U$.
\end{proposition}
\begin{proof}
Let $f_U=|U|^{-1}\sum_{u\in U}u$. It is an idempotent with image $S^U$, and it is self-adjoint for $[\,,\,]$ because $[us,s']=[s,u^{-1}s']$. So for $x\in S^U$ and $s\in S$ we have $[x,s]=[f_Ux,s]=[x,f_Us]$. Now let $x\in S^U$ be orthogonal to $S^U$, that is, in the radical of the restriction of $[\,,\,]$ to $S^U$. Then $[x,s]=[x,f_Us]=0$ for every $s\in S$, since $f_Us\in S^U$. So $x\in\mathrm{Rad}$. Conversely every $x\in S^U\cap\mathrm{Rad}$ is orthogonal to $S\supseteq S^U$. Hence the radical of the restricted form is $S^U\cap\mathrm{Rad}=\mathrm{Rad}^U$. By Proposition~\ref{propU} the $x_\tau$ form a basis of $S^U$, so the Gram rank is $\dim S^U-\dim\mathrm{Rad}^U=\dim(S/\mathrm{Rad})^U$, because taking $U$-fixed vectors preserves exact sequences ($|U|$ is invertible in $\KK$).
\end{proof}

By the Dipper--James theory, the $U$-fixed (equivalently $B$-fixed) vectors of the unipotent Specht module form the Specht module of the Hecke algebra $H_q(S_n)$ in characteristic $0$~\cite[(3.1)(vi), (2.15)(i), Example~5.2]{DJ89}. Proposition~\ref{propGramU} shows that the Gram ranks of the vectors $x_\tau$ are $\dim(D^\lambda)^U$, independently of the Hecke algebra. Computed values of these Gram ranks are given in Section~\ref{seccomp}, and their relation to the Hecke algebra is discussed in Section~\ref{secfurther}.

\section{Computations}\label{seccomp}

We computed $\dim D^\lambda$ for $\GL_n(2)$ for several partitions $\lambda$ and primes $\ell$. By Theorem~\ref{thmP}, the recursive family spans $S^\lambda_\KK$, so $\dim D^\lambda$ is the rank over $\FF_\ell$ of the Gram matrix $([b_i,b_j])$ of the recursive basis (Section~\ref{subD}). The vectors $b_i$ have integer entries, and the ranks are exact ranks of integer matrices modulo $\ell$. Table~\ref{tabD} shows the results.

\begin{table}[h]
\centering
\caption{$\dim D^\lambda$ for $\GL_n(2)$.}\label{tabD}
\begin{tabular}{c|c|cccc}
$\lambda$ & $\dim S^\lambda$ & $\ell=3$ & $5$ & $7$ & $31$\\\hline
$(2)$ & $2$ & $1$ & $2$ & $2$ & $2$\\
$(3)$ & $8$ & $7$ & $8$ & $3$ & $8$\\
$(2,1)$ & $6$ & $6$ & $6$ & $5$ & $6$\\
$(2,2)$ & $20$ & $7$ & $20$ & $19$ & $20$\\
$(2,1,1)$ & $14$ & $13$ & $13$ & $14$ & $14$\\
$(3,1)$ & $56$ & $35$ & $43$ & $56$ & $56$\\
$(4)$ & $64$ & $28$ & $21$ & $45$ & $64$\\
$(3,2)$ & $496$ & $217$ & $373$ & $495$ & $496$\\
$(3,1,1)$ & $280$ & $155$ & $280$ & $280$ & $251$\\
$(2,2,1)$ & $124$ & $124$ & $123$ & $94$ & $124$\\
$(4,1)$ & $960$ & $930$ & $960$ & $465$ & $709$\\
$(2,1,1,1)$ & $30$ & $30$ & $30$ & $30$ & $29$\\
\end{tabular}
\end{table}

When $\ell\nmid|\GL_n(2)|$ ($\ell=5,7,31$ for $n=2$; $\ell=5,31$ for $n=3$; $\ell=31$ for $n=4$), the Gram matrix is nondegenerate and $D^\lambda=S^\lambda$, as it must be. For $\GL_2(2)\cong S_3$ at $\ell=3$, the Steinberg module has a one-dimensional radical, so $\dim D^{(2)}=1=q-1$. For $q=3$ and $\ell=13$, which divides $|\GL_n(3)|$ for $n\ge3$, we found $\dim D^\lambda=16,\ 11,\ 89,\ 640$ for $\lambda=(3),(2,1),(2,2),(4)$, and full rank for the other shapes with $n\le4$.

The values of $\dim D^\lambda$ for $n\le4$ agree with the decomposition matrices in~\cite[pp.~106--108]{James84}, with one exception. James's labels are the transposes of ours, so his $D_\mu$ is our $D^{\mu'}$. For the two-column shapes (two-part partitions in James's labeling), the reducibility pattern also agrees with~\cite[Thm.~20.3]{James84}. The exception is the row $(1^4)$ of the matrix on p.~107, for $\ell\mid1+q$. For $q=2$ and $\ell=3$, that row gives, as we read it, $\dim D=64-1-7-35=21$ for the Steinberg module. We find $28$: the $\FF_3$-span of the translates of Steinberg's vector in the $315$ flags of $\FF_2^4$ has dimension $64$ and Gram rank $28$ (and $21$ and $45$ at $\ell=5,7$, as in the table). So for $q=2$ and $\ell=3$ our calculation indicates that the entry of that row in the column of $D_{(2,2)}$ should be $0$.

By Proposition~\ref{propGramU}, the Gram rank of the vectors $x_\tau$, $\tau$ a standard chain, is $\dim(D^\lambda)^U$. Table~\ref{tabU} gives these ranks for $\GL_n(2)$; a value of $0$ means that $D^\lambda$ has no $U$-fixed vectors.

\begin{table}[h]
\centering
\caption{$\dim(D^\lambda)^U$ for $\GL_n(2)$; here $f^\lambda$ is the number of standard tableaux of shape $\lambda$.}\label{tabU}
\begin{tabular}{c|c|cccc}
$\lambda$ & $f^\lambda$ & $\ell=3$ & $5$ & $7$ & $31$\\\hline
$(2,1)$ & $2$ & $2$ & $2$ & $1$ & $2$\\
$(3)$ & $1$ & $0$ & $1$ & $0$ & $1$\\
$(2,1,1)$ & $3$ & $2$ & $2$ & $3$ & $3$\\
$(3,1)$ & $3$ & $0$ & $1$ & $3$ & $3$\\
$(2,2)$ & $2$ & $0$ & $2$ & $1$ & $2$\\
$(4)$ & $1$ & $0$ & $0$ & $0$ & $1$\\
$(4,1)$ & $4$ & $0$ & $4$ & $0$ & $1$\\
$(3,2)$ & $5$ & $0$ & $1$ & $4$ & $5$\\
$(3,1,1)$ & $6$ & $0$ & $6$ & $6$ & $3$\\
$(2,2,1)$ & $5$ & $5$ & $4$ & $1$ & $5$\\
$(2,1,1,1)$ & $4$ & $4$ & $4$ & $4$ & $3$\\
\end{tabular}
\end{table}

These agree with the dimensions of the simple modules of the Hecke algebra $H_2(S_n)$ over $\FF_\ell$ that one gets from block theory, where $e=\mathrm{ord}_\ell(2)$ determines the blocks. For example, for $n=4$ and $\ell=7$ ($e=3$), the partitions $(3,1)$ and $(2,1,1)$ are $3$-cores, each giving a simple module of dimension $3$, and $(4),(2,2),(1^4)$ form a block with simple modules of dimensions $1$ and $1$. The table gives $3$, $3$, and $1$ for $(3,1)$, $(2,1,1)$, and $(2,2)$, and the trivial module $(1^4)$ accounts for the remaining $1$.

\section{Further directions and limitations}\label{secfurther}

The following questions remain open.

\begin{itemize}
\item \emph{A rule for $D^\lambda$.} Proposition~\ref{propD} reduces finding a basis of $D^\lambda$ to choosing the subfamilies $I_k$, that is, to describing the radical of the form layer by layer in the coordinates of the recursive basis. A selection rule needs more than the decomposition numbers, which give only the dimensions of the simple modules and the composition multiplicities of the Specht modules. For $S_n$ no uniform basis of $D^\lambda$ is known to us either (Section~\ref{subsym}).
\item \emph{The Hecke algebra.} In characteristic $0$, the $U$-fixed vectors of $S^\lambda$ form the Specht module of the Hecke algebra $H_q(S_n)$ (Section~\ref{subsym}). Proposition~\ref{propGramU} shows that the Gram rank of the vectors $x_\tau$ is $\dim(D^\lambda)^U$, and the examples of Section~\ref{seccomp} agree with the block theory of $H_2(S_n)$. We have not proved that the restriction of the form $[\,,\,]$ to the $U$-fixed vectors is the Gram form of the Hecke algebra Specht module, up to a unit. This would relate the Gram ranks of the $x_\tau$ to the simple modules of $H_q(S_n)$, and for $q\equiv1\pmod\ell$ to those of $S_n$. The fixed-vector Gram ranks of Section~\ref{seccomp} are computed for $q=2$, and the computation for $q=3$ and $\ell=13$ does not test this either. The symmetric group case needs $q\equiv1\pmod\ell$, for example $q=3$ and $\ell=2$.
\item \emph{Naturality.} The basis is adapted to the layers $\Omega_{Y_k}S^\lambda$ (the restriction to the mirabolic subgroup and the groups $\GL_{n-k}$) by construction (Lemma~\ref{lemlayers} and Proposition~\ref{propfilt}). The sets $\mathcal B_s$ of Proposition~\ref{propstd} regroup the same vectors by standardizing the chain. They are defined only combinatorially, and they depend on the choice of the matrices $h_w$ in the sets $R_k$ (Section~\ref{subnotation}). We do not know whether they are compatible with the module structure, in the way that the standard polytabloid basis of $S_n$ is compatible with restriction to $S_{n-1}$.
\item \emph{Bases of vectors $e_T$.} The recursive family is built from the generators $x_T$ by Fourier sums, and it need not consist of vectors $e_T$, although some of its vectors are of this form (for example the vectors $b_c$ of Example~\ref{exSt}). We do not give a combinatorial choice of tableaux $T$ for which the vectors $e_T$ form a basis of $S^\lambda$. This is the question stated at the beginning of the Introduction.
\item \emph{Other groups.} The finite orthogonal, symplectic, and unitary groups also have unipotent representations. Whether an analogous tableau construction and recursive basis exist for them is, to our knowledge, open.
\end{itemize}

\section*{Use of AI}
The construction, proofs, and computations in this paper were developed in extended collaboration with Claude (Anthropic), which proposed and checked arguments, wrote the code, and drafted the text. The author has independently verified the mathematics.

\end{document}